\documentclass[11pt]{article}
\usepackage{amssymb,amsmath,amsthm,mathtools}
\newcommand{\Mod}[1]{\ (\mathrm{mod}\ #1)}

\usepackage{graphicx,fullpage} 
\usepackage{url}

\usepackage[colorlinks,allcolors=blue]{hyperref}
\usepackage[numbers]{natbib}
\usepackage{longtable} \usepackage[ruled,linesnumbered]{algorithm2e}

\usepackage[shortlabels]{enumitem}
\setlist{
  listparindent=3em,
  itemindent=\parindent,
  parsep=0pt,
  leftmargin=\parindent
}

\newlist{nested}{enumerate}{5}

\setlist[nested]{
  nosep,    
  noitemsep,
  listparindent=2\parindent,
    leftmargin=2\parindent,
  parsep=0pt
}

\usepackage{caption}
\usepackage[dvipsnames]{xcolor}
\usepackage{soul}
\usepackage{verbatim}
\usepackage{comment}
\usepackage{changepage}

\newenvironment{indented}[1]
{   
    \begin{adjustwidth}{1cm}{1cm}
    #1
    \end{adjustwidth}
}

\usepackage{tikz-cd}

\usepackage{array}
\newcommand{\old}[1]{}
\mathchardef\mhyphen="2D 

\newtheorem{thm}{Theorem}[section]
\newtheorem{lemma}[thm]{Lemma}

\newtheorem{conj}[thm]{Conjecture}
\newtheorem{dfn}[thm]{Definition}

\newtheorem{conjecture}[thm]{Conjecture}

\newtheoremstyle{case}{}{}{}{}{}{:}{ }{}
\theoremstyle{case}
\newtheorem{case}{Case}
\newtheorem{obs}[case]{\bf{Observation}}

\usepackage{xspace} 

\newcommand{\f}{{\textit{PRF}}\xspace}
\newcommand{\fs}{{\textit{PRFs}}\xspace}

\newcommand{\Fq}{\mathbb{F}_q}

\DeclareMathAlphabet{\mathpzc}{OT1}{pzc}{m}{it}

\def\G{{\mathcal G}}

\def\P{{\mathcal P}}

\def\S{{\mathcal S}}
\def\T{{\mathcal T}}

\def\W{{\mathpzc W}}
\def\Y{{\mathpzc Y}}
\def\Z{{\mathpzc Z}}

\newcommand{\ie} [1] {\textit{i.e.,} #1}

\def\etal{{\it et~al.}\,}

\colorlet{Mycolor1}{green!50!}

\usepackage[dvipsnames]{xcolor}

\begin{document}

\title{Using  Permutation Rational Functions to Obtain Permutation Arrays with Large Hamming Distance.}

\author{Sergey Bereg\thanks{
Department of Computer Science,
University of Texas at Dallas,
Box 830688,
Richardson, TX 75083
USA. Research of the first author is supported in part by NSF award CCF-1718994. }
\and Brian Malouf$^*$
\and Linda Morales$^*$
\and Thomas Stanley$^*$
\and I. Hal Sudborough$^*$\thanks{Corresponding author: hal@utdallas.edu}
}

\providecommand{\keywords}[1]
{
  \small	
  \textbf{\textit{Keywords---}} #1
}
\maketitle

\begin{abstract}
    We consider permutation rational functions (PRFs), $V(x)/U(x)$, where both $V(x)$ and $U(x)$ are polynomials over a finite field $\Fq$. 
    Permutation rational functions have been the subject of several recent papers.
    Let $M(n,D)$ denote the maximum number of permutations on $n$ symbols with pairwise Hamming distance $D$. 
    Computing lower bounds for $M(n,D)$ is the subject of current research with applications in error correcting codes.
    Using \fs of specified degrees we  obtain improved lower bounds for $M(q,q-d)$ and $M(q+1,q-d)$, for prime powers $q$ and $d \in \{4,5,6,7,8,9\}$. 
\end{abstract}


\section{Introduction}
\label{sec:intro}
Permutation arrays (PAs) with large Hamming distance have applications in the design of error correcting codes. 
Lower bounds for the size of such permutation arrays are given, for example, in \cite{bms17,bls-18,bmmssw-19equiv,bmmms-19,bmms-19,bmmssPRFwaifi2020,chu2004,colbourn2004,gao13,jani18,jani15,Neri19,pavlidou2003power,wang17,Yang2006}.

Let $\pi$ and $\sigma$ be permutations over $n$ symbols.
The Hamming distance between $\pi$ and $\sigma$, denoted by $hd(\pi,\sigma)$, is the number of positions $x$ such that $\pi(x) \ne \sigma(x)$. 
Define the Hamming distance of a PA $A$, by $hd(A)=\min \{ hd(\pi,\sigma)\mid  \pi, \sigma \in A,~\pi \ne \sigma \}$. 
Let $M(n,D)$ denote the maximum number of permutations in any PA $A$ on $n$ symbols with Hamming distance $D$.

Let $\Fq$ denote the finite field with $q=p^m$ elements, where $p$ is prime and $m \ge 1$. 
The prime $p$ is called the {\em characteristic} of the field. 
A polynomial $V(x)$ over $\Fq$ is a {\em permutation  polynomial} 
if it permutes the elements of $\Fq$. 
Permutation polynomials have been studied for many decades, for example \cite{bmmssw-19equiv,Fan19,Fan19a,Fan19b,Hou15,li10,lidl93,Shallue-Wanless-pp-13}. 

This paper focuses on PAs derived from {\em permutation rational functions} (\fs), defined as follows:
\begin{dfn}
\label{dfn:fpp}
Let $V(x)$ and $U(x)$ be polynomials over $\Fq$, such that $\gcd(V(x),U(x))=1$.
Let  $\P^1(\Fq)$ denote $\Fq \cup \{\infty\}$.
If the rational function $V(x)/U(x)$ permutes the elements of  $\P^1(\Fq)$, then it is called a {\bf permutation rational function (\f)}. A \f has degree k \cite{Ding2020}, if the degree of the numerator (denominator) is $k$ and the degree of the denominator (numerator) is at most $k$.
\end{dfn}
Ferraguti and Micheli \cite{Ferraguti20} enumerated all \fs of degree 3. 
Ding and Zieve \cite{Ding2020} and Hou \cite{Hou2020} described all \fs of degree 4. 
Yang \etal \cite{Yang2006} were the first to use \fs to compute improved  lower bounds for $M(q,D)$ and $M(q+1,D)$ for some $q$ and $D$. 
In particular, they listed six lower bounds, namely, $M(18,14)$, $M(19,14)$, $M(19,15)$,  $M(19,16)$, $M(20,14)$, and $M(24,20)$.   
Of these, $M(19,14)\ge 65,322$ and $M(24,20)\ge 23,782$ are competitive. 
The rest are smaller than known lower bounds.
We have reproduced their bound for $M(19,14)$ experimentally. 
However, 
we have not been able to reproduce their result for $M(24,20)$. 
Instead we found $M(24,19) \ge 23,782$.

This paper extends our previous work on PAs with large Hamming distance \cite{bms17,bls-18,bmmssw-19equiv,bmmms-19,bmms-19} and \fs\cite{bmmssPRFwaifi2020} and gives improved lower bounds for $M(q,D)$ and $M(q+1,D)$ for various $D$. 
Our paper is organized as follows. 
In Section \ref{sec:counting-fpps} we summarize known information about the number of \fs of various degrees, and discuss techniques and new formulas for counting \fs for the purpose of Hamming distance.
In Section \ref{sec:hd-from-T_d(q)} we discuss Hamming distance properties of PAs on $\P^1(\Fq)$ derived from \fs, 
and give proofs of our new lower bounds for $M(q+1,q-d)$ and derive explicit formulas for $M(q+1,q-d)$, for $d=4,5,6$ and 7. 
In Section \ref{sec:hd-from-S_d(q)} we discuss Hamming distance properties of PAs on $\Fq$ created by contraction of permutations derived from \fs, 
and give proofs of our new lower bounds for $M(q,q-d)$ and derive explicit formulas for $q>7$ and $5 \le d \le 7$.
In Section \ref{sec:conjectures}, we briefly discuss our computational techniques and compute the number of \fs of degree 5, for $q \le 127$. 
These numbers are used for computing $M(q,q-d)$ and $M(q+1,q-d)$, for $d = 8,9$ and $q \le 127$.
We also give conjectured lower bounds for $M(q,q-8),~M(q,q-9),~M(q+1,q-8)$ and $M(q+1,q-9)$, for $q \ge 128$. In Section \ref{sec:conclusions}, we summarize our results and make additional conjectures.

\medskip
{\em Notation.} We use the following notation throughout this paper.
The letter $q$ represents a power of a prime and $\Fq$ denotes the finite field of $q$ elements, where $q=p^m$ for some prime $p$ and $m\ge 1$.
For notational clarity, we let $V,~U,~R$ and $S$ denote polynomials of degree $v,~u,~r$ and $s$, with coefficients $a_i, b_i, c_i$ and $d_i$ respectively,
That is, $V(x)=\sum_{i=0}^v{a_i x^i}$, $U(x)=\sum_{i=0}^u{b_i x^i}$,  $R(x)=\sum_{i=0}^r{c_i x^i}$, and $S(x)=\sum_{i=0}^s{d_i x^i}$, 
We denote \fs over $\P(\Fq)$ by the script symbols  $\W,~\Y$, and $\Z$. 
In other words $\W(x)=\frac{V(x)}{U(x)}= \sum_{i=0}^v{a_i x^i} / \sum_{i=0}^u{b_i x^i}$.

We use standard rules for evaluating $\W(\infty)$, which we list below for reference. 
In the last case, $a_v$ and $b_v$ are the high order coefficients of $V(x)$ and $U(x)$, respectively.
See \cite{bmmssPRFwaifi2020} for more detail.
Observe that when $v>u$, permutations generated by \fs over $\P^1(\Fq)$ can be viewed as permutations of $\Fq$ by eliminating $\infty$ from the domain.
  \begin{equation}
    \W(\infty)=
    \begin{cases}
        \infty,     &\mathrm{~when~} v>u, \label{eqn:W(infty)}\\
 
        0,          & \mathrm{~when~} v<u,\\
    
        a_v/b_v,    & \mathrm{~when~} v=u.
    \end{cases}
  \end{equation}


\section{Counting \fs}
\label{sec:counting-fpps}

In general, many of the  concepts and techniques discussed for polynomials over finite fields apply to \fs. 
A polynomial is {\it monic} if the coefficient of its high order term is 1.
A \f $\W(x) = V(x)/U(x)$ is {\it monic} if the coefficients of the high order terms of both $V(x)$ and $U(x)$ are 1.
Let $N_d(q)$ be the number of permutation polynomials of degree $d$ over $\Fq$ \cite{lidl93}. 
We generalize by defining $N_{v,u}(q)$, $N^m_{v,u}(q)$ and other useful concepts for \fs.
\begin{dfn} \label{dfn:counting}
\begin{enumerate}[(a)] 
\item \label{dfn:N_{vu}}
Let 
$\boldsymbol{P_{v,u}(q)}=\{V(x)/U(x) ~|~ V(x)$ has degree $v$ and $U(x)$ has degree $u$, and $V(x)/U(x)$ is a \f \}.
Let $\boldsymbol{\Pi_{v,u}(q)}$ be the set of permutations on $\P^1(\Fq)$ defined by the \fs in $P_{v,u}(q)$.
Let $\boldsymbol{N_{v,u}(q)}=|P_{v,u}(q)|=|\Pi_{v,u}(q)|$ be the number of \fs in $P_{v,u}(q)$. 
\item \label{dfn:N^m_{vu}}
Let $\boldsymbol{P^m_{v,u}(q)}$ be the set of monic \fs in $P_{v,u}(q)$, that is, the set of \fs where both $V(x)$ and $U(x)$ are monic.
Let  $\boldsymbol{\Pi^m_{v,u}(q)}$ be the set of permutations on $\P^1(\Fq)$ defined by the \fs in $P^m_{v,u}(q)$.
Let $\boldsymbol{N^m_{v,u}(q)}=|P^m_{v,u}(q)|=|\Pi^m_{v,u}(q)|$ be the number of monic \fs in $P_{v,u}(q)$.
\item \label{dfn:G_v(q)}
Let $\boldsymbol{\G_v(q)}=N_{v,v}(q)+ 2\sum_{u<v} N_{v,u}(q)$ be the number of degree $v$ \fs.
\end{enumerate}
\end{dfn}





\bigskip
As will soon be apparent, to obtain lower bounds for $M(n,D)$ we need to know the number of degree $k$ \fs of the form $\frac{\text{deg } k}{\text{deg } j}$ for a given $k$ and $j$.
For example, for all $q$, we need to know $N_{4,3}(q)$. 
Such specific values are not immediately evident from \cite{Ding2020}, \cite{Hou2020} or \cite{Hou20}.
We derive formulas and compute values of $N_{v,u}(q)$, for many additional values of $v,~u$ and $q$,
and use the computed values to give significantly improved lower bounds for $M(q,D)$ and $M(q+1,D)$ for various $D$. 

To begin, we show in Theorem \ref{thm:equal-num-denom} below that for any $v>1$, $N_{v,v}(q) = (q-1)\sum_{u < v} N_{v, u}(q)$.
$N_{v,v}(q)$ and $N_{v, u}(q)$ are used in the computation of $\S_d(q)$ and $\T_d(q)$ (defined below), and which are used in our lower bound theorems for $M(q,D)$ and $M(q+1,D)$.
In Theorem \ref{thm:G_v}, we provide a simpler formula for $\G_v(q)$ which we use in many proofs in the paper. 

\begin{lemma}
\label{lemma:invert-plus-c}
  If $V(x)/U(x)$ is a \f of degree $v$ over degree $u$ and $v > u$, then  $U(x)/V(x) + c$ for any $c \neq 0$, is a \f of degree $v$ over degree $v$.
\end{lemma}

\begin{proof}
  Assume $V(x)/U(x)$ is a \f of degree $v$ over degree $u$ and $v > u$. 
  Then $U(x)/V(x)$ is also a \f. We can express the addition of any constant $c$ to $U(x)/V(x)$ as a new \f $\frac{U(x) + cV(x)}{V(x)}$. 
  If $c$ is nonzero, then the coefficient of the $x^v$ term of $cV(x)$ will also be nonzero. 
  Since $v > u$, the resulting \f must be degree $v$ over degree $v$.
\end{proof}

\begin{thm}
\label{thm:equal-num-denom}
  For all $q$ and for any $v>1$, $N_{v,v}(q) = (q-1)\sum_{u < v} N_{v, u}(q)$.
\end{thm}

\begin{proof}
Let $v$ be arbitrary, and let $u<v$.
By Definition \ref{dfn:counting}\ref{dfn:N_{vu}} there are $N_{v, u}(q)$ \fs of degree $v$ over degree $u$ for any specific $u$.
By Lemma \ref{lemma:invert-plus-c}, each gives rise to a \f of degree $v$ over degree $v$ of the form $U(x) + cV(x)$. 
Since there are $q-1$ nonzero options for $c \in \Fq$, there are $(q-1)N_{v, u}(q)$ \fs of degree $v$ over degree $v$.
Thus $N_{v,v}(q)=(q-1)\sum_{u < v} N_{v, u}(q)$.

 To show each of these \fs is unique, assume we have two \fs $\frac{V(x)}{U(x)}$ and $\frac{R(x)}{S(x)}$ where $v > u$, $r > s$, and $v = r$, and $\frac{U(x)}{V(x)} + c = \frac{S(x)}{R(x)} + d$ for some constants $c$ and $d$. 
 We have
  \begin{align*}
    & \frac{U(x)}{V(x)} = \frac{S(x)}{R(x)} + c' \text{, where } c' = d - c \\
    & \frac{U(x)}{V(x)} = \frac{S(x) + c'R(x)}{R(x)} \\
    & U(x)R(x) = S(x)V(x) + c'R(x)V(x)
  \end{align*}
Note that the polynomial $U(x)R(x)$ has degree $u+r$ and the polynomial the polynomial $S(x)V(x) + c'R(x)V(x)$ has degree $r+v>u+r$ since $u < v = r$.
Thus equality holds only when $c' = 0$, in other words, when $c=d$.
That is, the \fs in the summation are unique.

To show that all \fs of degree $v$ over degree $v$ are accounted for, we show that any such \f can be transformed to a \f where the degree of the numerator is greater than the degree of the denominator by performing the operations in  reverse order. 
Consider an arbitrary \f $\frac{U(x)}{V(x)}$ such that $u = v$. 
Adding a constant $c$ gives us $\frac{U(x) + cV(x)}{V(x)}$. 
If we choose $c$ such that the leading coefficient of the numerator is $c = \frac{-b_u}{a_v}$ so that $b_u + ca_v = 0$, then the resulting \f $\frac{U'(x)}{V(x)}$ will have $u' < v$. 
Note that depending on the values of the other coefficients, it is possible that more than just the leading coefficient is zeroed out. 
If we then take the inverse \f $\frac{V(x)}{U'(x)}$, we will have a \f where the degree of the numerator is greater than the degree of the denominator, and thus $\frac{U(x)}{V(x)}$ can be found through Lemma \ref{lemma:invert-plus-c}.
\end{proof}

\begin{thm}
\label{thm:G_v}
For $q$ and for all $v>1$, $\G_v(q) = (q+1)\sum_{u<v} N_{v,u}(q)$.
\end{thm}

\begin{proof}
By Definition \ref{dfn:counting}\ref{dfn:G_v(q)}, $\G_v(q)=N_{v,v}(q)+ 2\sum_{u<v} N_{v,u}(q)$. 
Using Theorem \ref{thm:equal-num-denom} to substitute for $N_{v,v}(q)$ gives 
$\G_v(q)=(q-1)\sum_{u < v} N_{v, u}(q)+ 2\sum_{u<v} N_{v,u}(q) = (q+1)\sum_{u < v} N_{v, u}(q)$.
\end{proof}

\bigskip
We make the following straightforward observations, which are used in many of the proofs in this paper.

\begin{obs}    
\label{obs:N_d,0}
    $\boldsymbol{N_{d,0}(q)=N_d(q)}$. 
    \indented{Of specific interest in this paper are the formulas for  $N_{2,0}(q),~ N_{3,0}(q), ~N_{4,0}(q)$, and $N_{5,0}(q)$, the number of permutation polynomials of degree 2, 3, 4 and 5, respectively.
    We refer the reader to Table 2 of \cite{chu2004}. For the convenience of the reader, we reproduce these formulas  as  equations \ref{eqn:N_2,0} through \ref{eqn:N_5,0} below.}
\end{obs}

\begin{obs}    
\label{obs:N_vu=N_uv}
    $\boldsymbol{N_{u,v}(q)=N_{v,u}(q)}$.
    \indented{To see why, note that $\frac{V(x)}{U(x)}$ is a \f if and only if $\frac{U(x)}{V(x)}$ is also a \f. That is, if $(\alpha_0,\alpha_1, \dots ,\alpha_q)$ is a permutation of $\P^1(\Fq)$, then 
    $(\alpha_0^{-1},\alpha_1^{-1}, \dots ,\alpha_q^{-1})$ is also a permutation of $\P^1(\Fq)$.}
\end{obs}

\begin{obs}    
    \label{obs:N^m_vu} $\boldsymbol{N^m_{v,u}(q)=\frac{1}{q-1}N_{v,u}(q)}$.
    \indented{To see why, note that $\frac{V(x)}{U(x)}=\sum_{i=0}^v{a_i x^i} / U(x) =a_v x^v/U(x) + \sum_{i=0}^{v-1}{a_i x^i} / U(x)$.
    There are $q-1$ nonzero values for the high-order coefficient $a_v\in \Fq$, 
    so by the Pigeonhole Principle, there are $\frac{1}{q-1}N_{v,u}(q) $ monic \fs of degree $v$ over degree $u$.}
\end{obs}

\begin{obs}    
\label{obs:N_v1=0}
    $\boldsymbol{N_{v,1}(q) = 0 = N_{1,v}(q)},\mathbf{~for~all~} \boldsymbol{v>1}$. 
    \indented{To see why, consider any rational function  $\frac{V(x)}{U(x)}$ such that $u=1$.
    Note that $U(x)$ has one root, say $a$, that is, $U(a)=0$. 
    Then $\frac{V(a)}{U(a)}=\frac{V(a)}{0}=\infty$.
    Also since $v>u$,  $\frac{V(\infty)}{U(\infty)}=\infty$.
    Thus, $\infty$ appears in two positions, so $\frac{V(x)}{U(x)}$ is not a \f.
    Hence, $N_{v,1}(q) =0$, and so $N_{1,v}(q) =0$ as well.}
\end{obs}

\bigskip
We summarize known results 
for degree 1, 2, 3, 4 and degree 5 \fs. First, from \cite{chu2004} we have
\begin{align}
    N_{1,0}(q) &= q(q-1) \label{eqn:N_1,0} \mathrm{~~~~~~~for~all~} q,\\
    N_{2,0}(q) &= 
    \begin{cases}
    \label{eqn:N_2,0}
        q(q-1) &\mathrm{~when~} q=2^m, \\
        0  &  \mathrm{~otherwise.}
    \end{cases}\\ 
    N_{3,0}(q) &=
    \begin{cases}
    \label{eqn:N_3,0}
        \frac{1}{2}q(q-1)^2 + q(q-1) &\mathrm{~when~} q \equiv 0 \Mod 3, \\
        0  & \mathrm{~when~} q \equiv 1 \Mod 3, \\
        q^2(q-1)   & \mathrm{~when~} q \equiv 2 \Mod 3. 
    \end{cases}\\ 
    N_{4,0}(q) &= 
    \begin{cases}
    \label{eqn:N_4,0}
        0  &  \mathrm{~for~odd~} q>7,\\
        \frac{1}{3}q(q-1)(q^2+2) &\mathrm{~for~even~} q>7.
    \end{cases}\\ 
    N_{5,0}(q) &=
    \begin{cases} \label{eqn:N_5,0}
        \frac{1}{2}q^2(q-1)^2 + \frac{3}{4}q(q-1)^2+q(q-1), &\mathrm{~when~} q \equiv{0 \Mod 5}, \\     0, &\mathrm{~when~}  q \equiv{1 \Mod 5},\\
        q^3(q-1), &\mathrm{~when~} q \equiv{2,3 \Mod 5},\\
        q^2(q-1),
        &\mathrm{~when~} q \equiv{4 \Mod 5}.
    \end{cases}     
\end{align}    

\noindent For degree 3 \fs the following results are known \cite{Ferraguti20}:
\begin{align}
N_{3,2}(q) &= \tfrac{1}{2}q^2(q-1)^2, \mathrm{~~~~~~~for~all~} q, \label{eqn:N_3,2}\\ 
N_{3,3}(q) &=
        \begin{cases} \label{eqn:N_3,3}
        \frac{1}{2}(q^5-2q^4+2q^3-2q^2+q), &\mathrm{~when~} q \equiv{0 \Mod 3}, \\
        \frac{1}{2}q^2(q-1)^3, &\mathrm{~when~}  q \equiv{1 \Mod 3},\\
        \frac{1}{2}q^2(q-1)^2(q+1), &\mathrm{~when~} q \equiv{2 \Mod 3}.
        \end{cases}
\end{align}

\noindent For degree 4 \fs, it is known that \cite{Hou20} 
\begin{equation}
    N_{4,2}(q) = 0,  \mathrm{~for~all~} q > 7~ \label{eqn:N_4,2}
\end{equation}
As observed earlier,  $N_{v,u}(q) =N_{u,v}(q)$, so $N_{2,4}(q)= 0$ as well. 
Recall that $\G_k(q)$ is the number of \fs of degree $k$. The number of degree 4 \fs is \cite{Ding2020,Hou2020}
\begin{equation}
     \G_4(q) =
     \begin{cases} \label{eqn:G_4}
     \frac{1}{3}(q^3-q)^2
     , &\mathrm{~when~} q \mathrm{~is~odd~and~} q >7,\\
    \frac{1}{3}q(q-1)(q+2)(q^3+1)
    , &\mathrm{~when~} q \mathrm{~is~even~and~} q >7.\\
    \end{cases}
\end{equation}

\bigskip
We previously conjectured these same formulas for degree 3 and degree 4 \fs, for all $q$ \cite{bmmssPRFwaifi2020}. 
In addition, we previously conjectured the following lower bound for all $q \ge 19$: $N_{5,4}(q) \ge (q+1)q^3(q-1)^2/2$ \cite{bmmssPRFwaifi2020}. 
In Section \ref{sec:conjectures} of this paper we give new explicit formulas for $N_{5,4}(q)$ and $N_{5,5}(q)$, which we have experimentally verified for all $32 \le q \le 127$. and conjecture to be true, for all $q > 127$.

\bigskip
We now use Theorem \ref{thm:G_v}, Observation \ref{obs:N_v1=0} and equations \ref{eqn:N_1,0} through \ref{eqn:N_3,3} to derive expressions for $G_1(q),~G_2(q)$ and $G_3(q)$ which are used in several theorems in Section \ref{sec:hd-from-T_d(q)}.

\begin{align}
    \G_1(q) &= (q+1)N_{1,0}(q) = (q+1)q(q-1)= |PGL(2,q)|.\\    
    \G_2(q) &= (q+1)N_{2,0}(q) = 
        \begin{cases}
        (q+1)q(q-1) &\mathrm{~when~} q=2^m, \\
        0  &  \mathrm{~otherwise.}
    \end{cases}\\ 
    \G_3(q) &= (q+1)(N_{3,2}(q)+N_{3,0}(q)) \nonumber\\    
       &= \begin{cases}
       \label{eqn:G_3}
        (q+1)(\tfrac{1}{2}q^2(q-1)^2+\frac{1}{2}q(q-1)^2 + q(q-1) ), &\mathrm{~when~} q \equiv 0 \Mod 3, \\
        0,  & \mathrm{~when~} q \equiv 1 \Mod 3, \\
        (q+1)(\tfrac{1}{2}q^2(q-1)^2+q^2(q-1)),   & \mathrm{~when~} q \equiv 2 \Mod 3. 
    \end{cases}
\end{align}

\bigskip
\bigskip
In Theorem \ref{thm:N4,3} below, we derive a closed formula for $N_{4,3}(q)$ for all prime powers $q>7$.

\begin{thm}
\label{thm:N4,3}
For all prime powers $q>7$, $N_{4,3}(q) = \frac{1}{3} (q+1)q^2(q-1)^2$.
\end{thm}

\begin{proof}
By Theorem \ref{thm:G_v} and using Observation \ref{obs:N_v1=0} and equation \ref{eqn:N_4,2} to eliminate terms, the number of degree 4 \fs for $q>7$ is 
\begin{align*}
    \G_4(q) &= (q+1)(N_{4,3}(q)+N_{4,0}(q)).
\end{align*}
Solving for $N_{4,3}(q)$, we get
\[
N_{4,3}(q) = \tfrac{1}{q+1}\G_4(q)- N_{4,0}(q).
\]

For odd values of $q>7$,
Ding and Zieve \cite{Ding2020} proved that there are $\G(q)=\frac{1}{3}(q^3-q)^2=\frac{1}{3}(q+1)^2q^2(q-1)^2$ degree 4 PRFs. 
By equation \ref{eqn:N_4,0} we know that $N_{4,0}(q) = 0$ for odd $q$. 
Hence
\begin{align*}
    N_{4,3}(q) &=\tfrac{1}{q+1}\G(q) \\
    &= \tfrac{1}{3}(q+1)q^2(q-1)^2 
\end{align*} 

Now consider even values of $q>7$.
Ding and Zieve \cite{Ding2020} proved that there are $\G(q) = \tfrac{1}{3}q(q-1)(q+2)(q^3+1)$ degree 4 \fs.
By equation \ref{eqn:N_4,0}, $N_{4,0}(q) = \frac{1}{3}q(q-1)(q^2+2)$ when $q>7$ is even.
Hence 
\begin{align*}
  N_{4,3}(q) &=\tfrac{1}{3(q+1)}q(q-1)(q+2)(q^3+1)-\tfrac{1}{3}q(q-1)(q^2+2) \\
  &= \tfrac{1}{3}(q+1)q^2(q-1)^2.  
\end{align*}
\end{proof}

We define two new quantities,  $\T_d(q)$ and $\S_d(q)$, which yield better lower bounds for $M(q+1,q-d)$ and $M(q,q-d)$ than those given previously in \cite{bmmssPRFwaifi2020}.
Furthermore, by adding the term $N^m_{\frac{d+3}{2},\frac{d+1}{2}}(q)$ to $\S_d(q)$ and to $\T_d(q)$, we obtain new lower bounds for $M(q,q-d-1)$ and $M(q+1,q-d-1)$. 
The choices for $u$ and $v$ in the definitions are specifically tailored to meet certain degree requirements for the \fs, as will be explained in Sections \ref{sec:hd-from-T_d(q)} and \ref{sec:hd-from-S_d(q)}.

\begin{dfn}
\label{dfn:T_d(q)}
For odd integers d, let $\boldsymbol{\T_d(q)} =|\bigcup_{v,u} \Pi_{v,u}(q) |=\sum_{v,u}{N_{v,u}(q) }$,
for all $v,u \le (d+1)/2$. 
Recursively, we have
\begin{align}
\boldsymbol{\T_d(q)} =
    \begin{cases}         \label{eqn:T_d-recursive}
        \G_1(q)=N_{1,1}(q) +N_{1,0}(q) +N_{0,1}(q)  = |PGL(2,q)|, 
        & \mathrm{~when~} d=1,\\
        \T_{d-2}(q) + \G_{\frac{d+1}{2}}(q),
        & \mathrm{~when~} d \ge 3.
    \end{cases}
\end{align}

where $|PGL(2,q)| = (q+1)q(q-1)$.    
\end{dfn} \qed

\noindent For example, 
\begin{align}
\label{eqn:T_3}
    \T_3(q) &= 
    \begin{cases}
        2(q+1)q(q-1)=2q^3-2q &\mathrm{~when~} q=2^m, \\
        (q+1)q(q-1)=q^3-q &\mathrm{~otherwise}.
    \end{cases} 
\end{align}

\bigskip

\begin{dfn}
\label{dfn:S_d(q)}
    For odd integers d, let $\boldsymbol{\S_d(q)} = \sum_{v,u} N_{v,u}(q)$,
    where $v$ and $u$ are evaluated as 
  \begin{equation*}
    v,u=
    \begin{cases}
        v \le (d+1)/2,
        &\mathrm{~when~} v>u, \\
 
        u \le (d-3)/2, & \mathrm{~when~} v<u, \mathrm{~and}\\
        u,v \le (d-3)/2,& \mathrm{~when~} v=u.
    \end{cases}
  \end{equation*}
    
Recursively, we have
\begin{align}
\boldsymbol{\S_d(q)} =
    \begin{cases}                  N_{1,0}(q) = |AGL(1,q)| 
        & \mathrm{~when~} d=1,\\
        \S_{d-2}(q) + \sum_{u<k}N_{k,u}(q) + \sum_{v \le k-2}{N_{v,k-2}(q) }
        &\mathrm{~when~} d \ge 3,~ k=\frac{d+1}{2}.
    \end{cases}     \label{eqn:S_d-recursive}
\end{align}
where $|AGL(1,q)| = q(q-1)$. \qed
\end{dfn} 
\noindent 
For example, \begin{align} \label{eqn:S_3}
    \S_3(q) &= 
    \begin{cases}
        2q(q-1)=2q^2-2q &\mathrm{~when~} q=2^m,\\
        q(q-1)=q^2-q &\mathrm{~otherwise}.
    \end{cases}
\end{align}

\bigskip
\noindent We use these definitions in the proofs that follow to show for odd $d \ge 3$ and $q>d$  that: 
\begin{itemize}[left=12pt]
    \item $M(q+1,q-d) \ge \T_d(q)$ (Theorem \ref{thm:M(q+1,q-d)}), 
    \item $M(q+1,q-d-1) \ge \T_d(q)+N^m_{\frac{d+3}{2},\frac{d+1}{2}}(q)$ (Theorem \ref{thm:M(q+1,q-d-1)}),
    \item $M(q+1,q-4) \ge \T_3(q)+\max\{N^m_{3,2}(q),N_{3,0}(q)\}$ (Theorem \ref{thm:M(q+1,q-4)}),
    \item  $M(q,q-d) \ge \S_d(q)$ (Theorem \ref{thm:M(q,q-d)}), and
    \item $M(q,q-d-1) \ge \S_d(q)+N^m_{\frac{d+3}{2},\frac{d+1}{2}}(q)$ (Theorem \ref{thm:M(q,q-d-1)}). 
\end{itemize}
Using these results, we derive explicit formulas for $M(q+1,D)$ for $D=q-4,q-5,q-6$ and $q-7$ (see Theorems  \ref{thm:M(q+1,q-5)formulas},  \ref{thm:M(q+1,q-7)formulas},  \ref{thm:M(q+1,q-4)formulas} and
\ref{thm:M(q+1,q-6)formulas}),
and for $M(q,D)$ for $D=q-5,q-6$ and $q-7$ (see Theorems  \ref{thm:M(q,q-5)formulas},  \ref{thm:M(q,q-7)formulas} and
\ref{thm:M(q,q-6)formulas}).
In addition, we provide computational results to verify that 
$M(q+1,q-9) \ge \T_9$ (see Theorem \ref{thm:M(q+1,q-9)}), 
$M(q,q-9) \ge \S_9$ (see Theorem \ref{thm:M(q,q-9)}),
$M(q+1,q-8) \ge \T_7 + N^m_{5,4}$ (see Theorem \ref{thm:M(q+1,q-8)}), and 
$M(q,q-8) \ge \S_7 + N^m_{5,4}$ (see Theorem \ref{thm:M(q,q-8)}).


\section{New lower bounds for $M(q+1,q-D)$ }
\label{sec:hd-from-T_d(q)}

Recall that by Definition \ref{dfn:fpp}, $\gcd(V(x),U(x))=1$ for any \f. 
This property is implicit in our counting arguments for \fs.

We now discuss properties of \fs that are useful for improving lower bounds for $M(q+1,D)$ for various $D$. 
Some similar ideas were given in \cite{Yang2006}.
For the proofs in this section, we consider the \fs $\W(x)=\frac{V(x)}{U(x)}$ and $\Y(x)=\frac{R(x)}{S(x)}$ that permute the elements of $\P^1(\Fq)$ and such that $V(x)S(x) - U(x)R(x)$ is not a constant.
The degrees of the \fs $\W$ and $\Y$ need be not be the same. 

Note that the number of values $a \in \Fq$ such that $\W(a)=\Y(a)$ is given by the number of roots for the polynomial $V(x)S(x) - U(x)R(x)$.  
That is, the number of agreements between the permutations generated by the \fs $\W$ and $Y$ is given by the degree of the polynomial, which is  $\max\{v+s,u+r\}$.

\begin{dfn}
\begin{align*}
    \boldsymbol{\delta_1} & = \deg(V(x)S(x) - U(x)R(x)) \le \max\{v+s,u+r\}, \mathrm{~and}\\
    \boldsymbol{\delta_2} &=
        \begin{cases}
        1, &\mathrm{~when~} \W(\infty)=\Y(\infty), \\
        0, & \mathrm{~otherwise.}
        \end{cases}
  \end{align*}
\label{def:delta1,delta2}
\end{dfn}

\begin{lemma} \label{lemma:hd}
Let $\pi$ and $\sigma$ be the permutations of $\P^1(\Fq)$ generated by  $\W(x)$ and $\Y(x)$ respectively. Then 
for all $q,~hd(\pi,\sigma)\ge q+1-(\delta_1+\delta_2)$. 
\end{lemma}

\begin{proof}
It suffices to show that $\pi$ and $\sigma$ agree in at most $\delta_1+\delta_2$ positions.
If $a \in \Fq$ and $\W(a)=\Y(a)$ then $V(a)S(a) - U(a)R(a)=0$, that is, $a$ is a root.
There are at most $\delta_1$ roots of $V(x)S(x) - U(x)R(x)$.
Finally, $\pi$ and $\sigma$ agree at $\infty$ if and only if $\delta_2=1$.
\end{proof}




\noindent It follows also that permutations corresponding to different \fs are different, because the permutations have non-trivial Hamming distance.

\medskip

Recall that $\Pi_{v,u}(q)$ is the set of permutations over $\P^1(\Fq)$ defined by \fs of the form $\W(x)=\frac{V(x)}{U(x)}$.
In  Lemma \ref{lemma:hd(Pi)}, we prove that the set of permutations $\bigcup_{v,u} \Pi_{v,u}(q)$, for $v,u \leq (d+1)/2$,
has Hamming distance at least $q-d$.
We use this in
Theorem \ref{thm:M(q+1,q-d)} to show that $M(q+1,q-d) \geq \T_d(q)$.
(This is a simplified version of an earlier proof in \cite{bmmssPRFwaifi2020}.)

\begin{lemma}
\label{lemma:hd(Pi)}
For $q>7$ and for odd $d \ge 3$, let $u,v \leq (d+1)/2$.
Then $hd(\bigcup_{v,u} \Pi_{v,u}(q)) \ge q-d.$
\end{lemma}

\begin{proof}
Let $\pi$ and $\sigma$ be distinct permutations in $\bigcup_{v,u} \Pi_{v,u}(q)$, for $u,v \leq (d+1)/2$. 
Let $\W(x)=\frac{V(x)}{U(x)}$ and $\Y(x)=\frac{R(x)}{S(x)}$ be the \fs that generate $\pi$ and $\sigma$, respectively. 
By Lemma \ref{lemma:hd}, it suffices to show that $\delta_1+\delta_2\le d+1$.
Note that $\delta_1=\deg(V(x)S(x) - U(x)R(x))\le d+1$. 
The lemma follows immediately if $\delta_2=0$,
so suppose instead that $\delta_2=1$, \ie $\W(\infty)=\Y(\infty)$.
It suffices to show that in this case, $\delta_1\le d$.

\begin{indented}

\noindent{\bf Case 1}. $\W(\infty)=\Y(\infty)=\infty$. 

By equation \ref{eqn:W(infty)}, $\W(\infty)=\Y(\infty)=\infty$ implies that $v>u$ and $r>s$.
So by Definition \ref{def:delta1,delta2}, $\delta_1 \le \max \{(\frac{d+1}{2}+\frac{d-1}{2}),(\frac{d-1}{2}+\frac{d+1}{2}) \} = d$.

\noindent{\bf Case 2}.  
$\W(\infty)=\Y(\infty)=0$. 

In this case, $v<u$ and $r<s$, so $\delta_1\le d$.

\noindent{\bf Case 3}.  
$\W(\infty)=\Y(\infty)\notin \{0,\infty\}$. 

In this case, $v=u, r=s$ and $a_v/b_u=c_r/d_s$. 
Thus the degree of the $v+s=u+r$ term in the polynomial $V(x)S(x) - U(x)R(x)$, namely $(a_vd_s-b_uc_r)x^{v+s}$, has a coefficient equal to zero, ensuring that the degree of the polynomial is less than $v+s$.
So $\delta_1 < v+s\le d+1$, \ie{$\delta_1 \le d$}.
\end{indented}
\end{proof}

\begin{thm}
\label{thm:M(q+1,q-d)}
For odd $d \ge 3$ and $q \ge d+2$, $M(q+1,q-d) \ge \T_d(q).$
\end{thm}

\begin{proof}
This follows from Lemma \ref{lemma:hd(Pi)}, because by definition, $\T_d(q) =|\bigcup_{v,u} \Pi_{v,u}(q)|$ for $u,v \leq (d+1)/2$.
\end{proof}

For example, for $q=47$ and $d=7$, Theorem \ref{thm:M(q+1,q-d)} gives $M(48,40) \ge
3,781,770,400$, which is an improved lower bound. The previous lower bound of
9,655,492 is obtained from permutation polynomials for $q=47$ \cite{bmmssw-19equiv}.
For q = 47 and d = 9, Theorem \ref{thm:M(q+1,q-d)} gives $M(48,38) \ge 262,478,007,808$, which is an improved lower bound, where the previous lower bound was 21,442,716, using permutation polynomials \cite{bmmssw-19equiv}  

\bigskip
In Theorems \ref{thm:M(q+1,q-5)formulas} and \ref{thm:M(q+1,q-7)formulas} below, we derive explicit formulas for $\T_5(q)$ and $\T_7(q)$ for $q>7$, by substituting known formulas for each term in the definitions of $\T_5(q)$ and $\T_7(q)$. 
These formulas are based on the mod 3 congruence classes of $q$.


\begin{thm}
\label{thm:M(q+1,q-5)formulas}
For $q>9$
    \begin{equation*}
    M(q+1,q-5) \ge \T_5(q)=
    \begin{cases}
        \frac{1}{2}(q^5+2q^3-3q),& \mathrm{~for~} q \equiv{0 \Mod 3},\\
        \frac{1}{2}(q^5-q^4+q^3+q^2-2q), & \mathrm{~for~odd~} q \equiv{1\Mod 3}, \\
        \frac{1}{2}(q^5-q^4+3q^3+q^2-4q), & \mathrm{~for~even~} q \equiv{1 \Mod 3},\\
        \frac{1}{2}(q^5+q^4+q^3-q^2-2q), & \mathrm{~for~odd~} q \equiv{2 \Mod 3},\\
        \frac{1}{2}(q^5+q^4+3q^3-q^2-4q), & \mathrm{~for~even~}q \equiv{2 \Mod 3}.
    \end{cases}
\end{equation*}
\end{thm}

\begin{proof}
We use Theorem \ref{thm:M(q+1,q-d)}.
We provide a detailed example of the derivation of the formula for $\T_5(q)$ when $q$ is odd and  $q \equiv{2 \Mod 3}$.
\begin{align*}
    \T_5(q) &= \T_3 + \G_3 &\mathrm{(Equation~14)}\\
    &= (q+1)q(q-1)+ (q+1)(\tfrac{1}{2}q^2(q-1)^2+q^2(q-1))
    &\mathrm{(Equations~ 13 ~and ~15)}\\ 
    &= \tfrac{1}{2}(q^5+q^4+q^3-q^2-2q).
\end{align*} 
Derivation of the other formulas is similar, using terms appropriate for each congruence class. 
\end{proof}

\begin{thm}
\label{thm:M(q+1,q-7)formulas}
For $q>9$
    \begin{equation*}
    M(q+1,q-7) \ge \T_7(q) = 
    \begin{cases}
        \frac{1}{6}(2q^6+3q^5-4q^4+6q^3+2q^2-9q),& \mathrm{~for~} q \equiv{0 \Mod 3},\\
        \frac{1}{6}(2q^6+3q^5-7q^4+3q^3+5q^2-6q), & \mathrm{~for~odd~} q \equiv{1\Mod 3}, \\
        \frac{1}{6}(2q^6+5q^5-7q^4+11q^3+5q^2-16q), & \mathrm{~for~even~} q \equiv{1 \Mod 3},\\
        \frac{1}{6} (2q^6+3q^5-q^4+3q^3-q^2-6q), & \mathrm{~for~odd~} q \equiv{2 \Mod 3},\\
        \frac{1}{6}(2q^6+5q^5-q^4+11q^3-q^2-16q), & \mathrm{~for~even~}q \equiv{2 \Mod 3}.
    \end{cases}
\end{equation*}
\end{thm}

\begin{proof}
We use Theorems \ref{thm:M(q+1,q-d)} and \ref{thm:M(q+1,q-5)formulas}.
We provide a detailed example of the derivation of the formula for $\T_7(q)$ when $q$ is even and  $q \equiv{2 \Mod 3}$.

\begin{align*}
    \T_7(q) &=\T_5(q)+\G_4(q)&\mathrm{(Equation ~14)}\\
    &= \frac{1}{2}(q^5+q^4+3q^3-q^2-4q) + \frac{1}{3}q(q-1)(q+2)(q^3+1)\\
&&\mathrm{(Theorem~3.5,~Equation~10)}\\
    &=\tfrac{1}{6} (2q^6+5q^5-q^4+11q^3-q^2-16q).
\end{align*}

Derivation of the other formulas is similar, using terms appropriate for each congruence class. \end{proof}


Theorems \ref{thm:M(q+1,q-d)}, \ref{thm:M(q+1,q-5)formulas}, and \ref{thm:M(q+1,q-7)formulas} are appropriate when the desired Hamming distance is $q-d$ is where $d$ is odd.
If we relax the Hamming distance constraint, we can consider an additional set of \fs, to obtain results for $M(q+1,q-d-1)$, \ie{when the desired Hamming distance is $q-d-1$ and, hence,  $d-1$ is even}.
By Definition \ref{dfn:counting}\ref{dfn:N^m_{vu}}, $P^m_{r,s}(q)$ is the set of monic \fs in $P_{r,s}(q)$, 
and the size of this set is $N^m_{r,s}(q).$
Let  $v,u,s \le \frac{d+1}{2}$ and $r=\frac{d+3}{2}$.
Lemma \ref{lemma:hd(monicPi)} shows that $\Pi^m_{r,s}(q)$, the set of permutations  generated by the monic \fs 
$P^m_{r,s}(q)$ 
has Hamming distance $q-d-1$.
Lemma \ref{lemma:hd(Pi,monicPi)} shows that permutations in the union of the two sets $\bigcup_{v,u} \Pi_{v,u}(q)$  and $\Pi^m_{r,s}(q)$ have pairwise Hamming distance at least $q-d-1$.
Theorem \ref{thm:M(q+1,q-d-1)} uses these results to give the new lower bound $M(q+1,q-d-1) \ge \T_d(q)+N^m_{\frac{d+3}{2},\frac{d+1}{2}}(q)$, for $q>7$ and for $d \ge 3$.
A special case is presented by $d=3$, and in Theorem \ref{thm:M(q+1,q-4)} we show that $M(q+1,q-4)$ has a better lower bound than Theorem \ref{thm:M(q+1,q-d-1)} would suggest. 
We give explicit formulas for $M(q+1,q-4)$ in Theorem \ref{thm:M(q+1,q-4)formulas}, and for $M(q+1,q-6)$ in Theorem \ref{thm:M(q+1,q-6)formulas}.

\begin{lemma}
\label{lemma:hd(monicPi)}
For odd $d \ge 3$ and for $q \ge d+3$, let $r=(d+3)/2$ and $s=(d+1)/2$.
Then $hd(\Pi^m_{r,s}(q)) \ge q-d-1.$
\end{lemma}

\begin{proof}
Let $\W(x)=\frac{V(x)}{U(x)}$ and $\Y(x)=\frac{R(x)}{S(x)}$ be distinct \fs in $P^m_{r,s}(q)$.
By Lemma \ref{lemma:hd}, it suffices to show that $\delta_1+\delta_2\le d+2$.
Note that $\delta_1=\deg(V(x)S(x) - U(x)R(x))\le d+2$. 
However, since the \fs $P^m_{r,s}(q)$ are monic, the high order term of the polynomial
$(V(x)S(x) - U(x)R(x))$ is zero, so the in fact, $\delta_1=\deg(V(x)S(x) - U(x)R(x))\le d+1$.
Note also that since $r>s, ~\W(\infty)=\Y(\infty)=\infty$.
Thus $\delta_2=1$, and the lemma follows.
\end{proof}

\begin{lemma}
\label{lemma:hd(Pi,monicPi)}
For odd $d \ge 3$ and for $q \ge d+3$, let $u,v,s \leq (d+1)/2$, and let $r \leq (d+3)/2$.
Then $hd(\bigcup_{v,u} \Pi_{v,u}(q),\Pi^m_{r,s}(q)) \ge q-d-1.$
\end{lemma}

\begin{proof}
Let $\pi \in \bigcup_{v,u} \Pi_{v,u}(q)$  and $\sigma \in \Pi^m_{r,s}(q)$.
Let $\W(x)=\frac{V(x)}{U(x)}$ and $\Y(x)=\frac{R(x)}{S(x)}$ be the \fs that generate $\pi$ and $\sigma$, respectively.
We show that $hd(\pi,\sigma) \geq q-d-1.$
By Lemma \ref{lemma:hd} it suffices to show that $\delta_1 + \delta_2 \leq d+2$.
Note that $\delta_1 \leq d+2$.
The theorem follows immediately if $\delta_2=0$.
So suppose $\delta_2=1$, \ie{$\W(\infty)=\Y(\infty)$}.
We show that $\delta_1 \leq d+1$.

\begin{indented}

\noindent{\bf Case 1}. $\W(\infty)=\Y(\infty)=\infty$. 

In this case, $v>u$ and $r>s$, so $\delta_1\le \frac{d+3}{2} + \frac{d-1}{2}=d+1$.

\noindent{\bf Case 2}.  
$\W(\infty)=\Y(\infty)=0$. 

Then $v<u$ and $r<s$, 
so $\delta_1\le \frac{d+1}{2} + \frac{d-1}{2}=d < d+1$.

\noindent{\bf Case 3}.  
$\W(\infty)=\Y(\infty)\notin \{0,\infty\}$. 

In this case, $v=u, r=s$ and $a_v/b_u=c_r/d_s$, so $\delta_1 \le \frac{d+1}{2} + \frac{d+1}{2}=d+1$.
\end{indented}
\end{proof}

\begin{thm}
\label{thm:M(q+1,q-d-1)}
For $q>7$ and odd $d \ge 5$, $M(q+1,q-d-1) \ge \T_d(q)+N^m_{\frac{d+3}{2},\frac{d+1}{2}}(q).$
\end{thm}

\begin{proof}
This follows from Lemmas \ref{lemma:hd(Pi)},~\ref{lemma:hd(monicPi)} and \ref{lemma:hd(Pi,monicPi)}, since $\T_d(q)=|\bigcup_{v,u} \Pi_{v,u}(q)|$, and $N^m_{r,s}(q)=|\Pi^m_{r,s}(q)|$, for $u,v,s \leq (d+1)/2$, and $r \leq (d+3)/2$.
\end{proof}

\bigskip
Theorem \ref{thm:M(q+1,q-d-1)} applies specifically for $d \ge 5$, where $N^m_{\frac{d+3}{2},\frac{d+1}{2}}(q)$ is consistently larger than $N_{\frac{d+3}{2},0}(q)$. 
However for $d=3$, this is not the case for $q>7$.
We use this to fine tune the upper bound for $M(q+1,q-4)$ in Lemmas \ref{lemma:hd(Pi,deg3PPs)}, \ref{thm:M(q+1,q-4)} and Theorem \ref{thm:M(q+1,q-4)formulas} below.
Note that when $d=3,~ u,v,s \le \frac{d+1}{2}=2$ and $r \le \frac{d+1}{2}=3$. 

\begin{lemma}
\label{lemma:hd(Pi,deg3PPs)}
Let $u,v \leq 2$.
Then for all $q,~hd(\bigcup_{v,u} \Pi_{v,u}(q),\Pi_{3,0}(q)) \ge q-4.$
\end{lemma}

\begin{proof}
We use the fact that  $hd(\Pi_{3,0}(q)) \ge q-3$ \cite{chu2004} .
Let $\pi \in \bigcup_{v,u} \Pi_{v,u}(q)$  and $\sigma \in \Pi_{3,0}(q)$.
Let $\W(x)=\frac{V(x)}{U(x)}$ and $\Y(x)=\frac{R(x)}{S(x)}$ 
be the \fs that generate $\pi,~  \sigma$, 
respectively.
We show that $hd(\pi,\sigma) \geq q-4.$
By Lemma \ref{lemma:hd} it suffices to show that $\delta_1 + \delta_2 \leq 5$.
First  note that $\delta_1 \le \max\{v+2,u+r\}=\max\{2+0,2+3\} \leq 5$.
The theorem follows immediately if $\delta_2=0$.
So suppose $\delta_2=1$, \ie{$\W(\infty)=\Y(\infty)$}.
Since $r=3$ and $s=0$, $r$ is always greater than $s$ so $\Y(\infty)=\infty$. 
Hence $\W(\infty)=\Y(\infty)=\infty$ implies that $v>u$.
This means that $u \leq\ 1$, so $\delta_1=\max\{2+0,1+3\} \leq 4$. The lemma follows.
\end{proof}

\begin{thm}
\label{thm:M(q+1,q-4)}
For $q>7,~M(q+1,q-4) \ge \T_3(q)+\max\{N^m_{3,2}(q),N_{3,0}(q)\}.$
\end{thm}

\begin{proof}
We use Lemmas \ref{lemma:hd(Pi,monicPi)} and \ref{lemma:hd(Pi,deg3PPs)}.
First we note that by Lemma \ref{lemma:hd(Pi,monicPi)}, $hd(\bigcup_{v,u \le 2} \Pi_{v,u}(q),\Pi^m_{3,2}(q)) \ge q-4$.
To obtain the largest set of permutations on $q+1$ symbols with  Hamming distance $q-4$, we add to $\T_3(q)$ either $N_{3,0}(q)$ or $N^m_{3,2}(q)$, whichever is larger. 
For $d=3$, we use equation \ref{eqn:N_3,2} and Observation \ref{obs:N^m_vu} to obtain $N^m_{\frac{d+3}{2},\frac{d+1}{2}}(q)= N^m_{3,2}(q)=(\frac{1}{q-1})N_{3,2}(q)=\frac{1}{2}q^2(q-1)$.

When $q \equiv{1 \Mod 3},~ N_{3,0}(q)=0$, so clearly in this case, $N^m_{3,2}(q) > N_{3,0}(q)$.
However, when $q \equiv 0,2 \Mod 3$, a better result can be achieved by using $N_{3,0}(q)$ instead of $N^m_{3,2}(q)$.
In particular, from equation \ref{eqn:N_3,0} when $q \equiv{0 \Mod 3},~ N_{3,0}(q)=\frac{1}{2}q(q-1)^2 + q(q-1)$,
and when $q \equiv{2 \Mod 3},~N_{3,0}(q)=q^2(q-1)$. 
Thus in both cases, $N^m_{3,2}(q) < N_{3,0}(q)$

The theorem follows from Lemma \ref{lemma:hd(Pi,monicPi)}, and Lemma \ref{lemma:hd(Pi,deg3PPs)} since $\T_3(q)=|\bigcup_{v,u} \Pi_{v,u}(q)|$ for $u,v \leq 2$, $N^m_{3,2}(q)=|\Pi^m_{3,2}(q)|$,  and $N_{3,0}(q)=|\Pi_{3,0}(q)|$.
\end{proof}

\begin{thm}
\label{thm:M(q+1,q-4)formulas}
    For $q>7$,
    \begin{equation*}
    M(q+1,q-4) \ge \T_3(q)+\max\{N^m_{3,2}(q),N_{3,0}(q)\}=
    \begin{cases}
        \frac{3}{2}(q^3-q),& \mathrm{~for~} q \equiv{0 \Mod 3},\\
        \frac{1}{2} (3q^3-q^2-2q), & \mathrm{~for~odd~} q \equiv{1\Mod 3}, \\
        \frac{1}{2}(5q^3-q^2-4q), & \mathrm{~for~even~} q \equiv{1 \Mod 3},\\
        2q^3-q^2-q, & \mathrm{~for~odd~} q \equiv{2 \Mod 3},\\
        3q^3-q^2-2q, & \mathrm{~for~even~}q \equiv{2 \Mod 3}.
    \end{cases}
\end{equation*}
\end{thm}

\begin{proof}
We use Theorem \ref{thm:M(q+1,q-4)} and the formulas for $T_3(q),~N_{3,0}(q)$ and $N^m_{3,2}(q)$ from equations \ref{eqn:T_3}, \ref{eqn:N_3,0}, and \ref{eqn:N_3,2} along with Observation \ref{obs:N^m_vu}.
The results are obtained by adding appropriate terms for each congruence class. 
By Theorem \ref{thm:M(q+1,q-4)}, the Hamming distance is at least $q-4$. \end{proof}


For example, for $q=47$, Theorem \ref{thm:M(q+1,q-4)formulas} gives $M(48,43) \ge
205,390$, which is an improved lower bound. 

\begin{thm}
\label{thm:M(q+1,q-6)formulas}
    For $q>9$, 
    \begin{equation*}
    M(q+1,q-6) \ge \T_5(q)+N^m_{4,3}(q) =
    \begin{cases}
        \frac{1}{6}(3q^5+2q^4+6q^3-2q^2-9q),& \mathrm{~for~} q \equiv{0 \Mod 3},\\
        \frac{1}{6}(3q^5-q^4+3q^3+q^2-6q), & \mathrm{~for~odd~} q \equiv{1\Mod 3}, \\
        \frac{1}{6}(3q^5-q^4+9q^3+q^2-12q), & \mathrm{~for~even~} q \equiv{1 \Mod 3},\\
        \frac{1}{6} (3q^5+5q^4+3q^3-5q^2-6q), & \mathrm{~for~odd~} q \equiv{2 \Mod 3},\\
        \frac{1}{6} (3q^5+5q^4+9q^3-5q^2-12q), & \mathrm{~for~even~}q \equiv{2 \Mod 3}.
    \end{cases}
\end{equation*}
\end{thm}

\begin{proof}
We use Theorem \ref{thm:M(q+1,q-d-1)}.
The formulas for $\T_5(q)$ based on the (mod 3) congruence classes of $q$ are given by Theorem \ref{thm:M(q+1,q-5)formulas}.
The formula for $N^m_{4,3}(q)$ is given by Theorem \ref{thm:N4,3} and Observation \ref{obs:N^m_vu}.
\end{proof}

For example, for $q=47$, Theorem \ref{thm:M(q+1,q-6)formulas} gives $M(48,41) \ge
118,788,928$, which is an improved lower bound. 


\section{New lower bounds for $M(q,q-D)$ }
\label{sec:hd-from-S_d(q)}

We now consider permutations on $\Fq$ (not $\P^1(\Fq)$) in order to derive new lower bounds for $M(q,q-D)$.
Such permutations can be created from permutations on $\P^1(\Fq)$ by an operation called {\em contraction} \cite{bls-18} which eliminates occurrences of the symbol $\infty$ in permutations on $\P^1(\Fq)$.
For the reader's convenience, we give a definition of the contraction operation here.
A more detailed discussion of contraction is given in \cite{bls-18} and \cite{bmmms-19}.
In particular, we note that while the operation of contraction changes both the length of a permutation and the Hamming distance of a set of permutations, the number of permutations in the set is unchanged.

In this section we let $\pi$ and $\sigma$ be permutations on $\P^1(\Fq)$ generated by the \fs $\W(x)$ = $\frac{V(x)}{U(x)}$ and
$\Y(x) = \frac{R(x)}{S(x)}$, respectively.
Let $\pi'$ and $\sigma'$ denote the  permutations on $\Fq$ generated by the operation of contraction on $\pi$ and $\sigma$, and let 
$\bigcup_{v,u} \Pi'_{v,u}(q)$ and 
$\Pi^{'m}_{u,v}(q)$ denote the PAs on $\Fq$ generated by the operation of contraction on the PAs 
$\bigcup_{v,u} \Pi_{v,u}(q)$ and $\Pi^m_{u,v}(q)$, respectively.


\begin{dfn}
\label{def:contraction}
Let $\pi$ be any permutation on $\P^1(\Fq)$.
The {\bf contraction} operation converts $\pi$ to a new permutation $\pi'$ on $\Fq$ by the following rules:
\begin{enumerate}
    \item If $\pi(\infty) = \infty$, then  simply eliminate the symbol $\infty$ creating a permutation $\pi'$  on $\Fq$.
    \item If $\pi(\infty)= a$, with $a \in \Fq$, then exchange the symbol $\infty$ wherever it occurs in $\pi$ with the symbol $a$. 
    This moves the symbol $\infty$ to the last position in the permutation, so it can be eliminated, creating a permutation $\pi'$ on $\Fq$.
\end{enumerate}    
\noindent Define $\boldsymbol{\delta_3}$ to be the number of new agreements created by the operation of contraction on two permutations $\pi$ and $\sigma$ on $\P^1(\Fq)$.
\end{dfn}

\begin{lemma}
\label{lemma:delta3}
\begin{equation*}
    \boldsymbol{\delta_3} \leq 
        \begin{cases}
            0 & \mathrm{when~} \pi(\infty)=\sigma(\infty) = \infty \\
            1 & \mathrm{when~} \pi(\infty) = \infty, 
            \mathrm{~and~} \sigma(\infty) = a \mathrm{~for~some~} a \in \Fq
            \\
            2 & \mathrm{when~} \pi(\infty) = b \mathrm{~for~some~} b, \in \Fq \mathrm{~and~} \sigma(\infty) = a \mathrm{~for~some~} a \in \Fq \\
        \end{cases}.
\end{equation*}
\end{lemma}

\begin{proof}
Let $\pi'$ and $\sigma'$
be two  permutations on $\Fq$, created by the contraction two distinct permutations $\pi$ and $\sigma$ on $\P^1(\Fq)$.
We consider the number of new agreements created by contraction.
\begin{enumerate} [{Case} 1{:}, leftmargin=0.5in]
    \item If $\pi'$ and $\sigma'$ are both created by Rule 1 of Definition \ref{def:contraction}, then no new agreements are created, so $\delta_3 = 0$.
    \item Suppose that $\pi'$ is created by Rule 1 and $\sigma'$ is created by Rule 2 (or vice-versa).
    Suppose also that $\sigma(x)=\infty$ for some $x \in \Fq$.
    By Rule 2, $\sigma'(x)=a$ which would result in  one new  agreement if and only if $\pi(x)=a$.
    Thus, this case, $\delta_3 \leq 1$.
    \item If $\pi'$ and $\sigma'$ are both created by Rule 2, then at most two new agreements could result.
    This worst case happens if and only if $\sigma(x)=\infty$, $\sigma(y)=b$, $\sigma(\infty)=a$,  $\pi(y)=\infty$,
    $\pi(x)=a$ and $\pi(\infty)=b$  for some $x,y,a,b \in \Fq, ~x \neq y$.
    So in this case, $\delta_3 \leq 2$.
\end{enumerate} \end{proof}

\begin{lemma}
\label{lemma:hd-contraction}
$hd(\pi',\sigma') \geq q-(\delta_1 + \delta_3)$.
\end{lemma}

\begin{proof}
It suffices to show that $\pi'$ and $\sigma'$ agree in at most $\delta_1+\delta_3$ positions.
If $x\in \Fq$ and $\W(x)=\Y(x)$ then $V(x)S(x) - U(x)R(x)=0$, that is, $x$ is a root.
There are at most $\delta_1$ roots of $V(x)S(x) - U(x)R(x)$, that is, there are at most $\delta_1$ positions where $\pi$ and $\sigma$ agree (before contraction).
Finally, by Lemma \ref{lemma:delta3}, contraction creates at most $\delta_3$ new agreements.
\end{proof}

In  Lemma \ref{lemma:hd(Pi-contraction)}, we prove that the set of permutations $\bigcup_{v,u} \Pi'_{v,u}(q)$ for $v,u$ as given, has Hamming distance at least $q-d$.
We use this in Theorem \ref{thm:M(q,q-d)} to show that $M(q,q-d) \geq \S_d(q)$.
(This is a simplified version of an earlier proof in \cite{bmmssPRFwaifi2020}.)

\begin{lemma}
\label{lemma:hd(Pi-contraction)}
For odd $d \ge 3$, let $v$ and $u$ be defined by 
  \begin{equation*}
    v,u=
    \begin{cases}
        v \le (d+1)/2,
        &\mathrm{~when~} v>u, \\
 
        u \le (d-3)/2, & \mathrm{~when~} v<u, \mathrm{~and}\\
        u,v \le (d-3)/2,& \mathrm{~when~} v=u.
    \end{cases}
  \end{equation*} 

 Then $\bigcup_{v,u} \Pi'_{v,u}(q)  \ge q-d.$
\end{lemma}

\begin{proof}
Let $\W(x)$ = $\frac{V(x)}{U(x)}$ and
$\Y(x) = \frac{R(x)}{S(x)}$ be distinct \fs which generate permutations in $\bigcup_{v,u} \Pi_{v,u}(q)$. 
By Lemma \ref{lemma:hd-contraction}, it suffices to show that $\delta_1+\delta_3\le d$.
We do a proof by cases based on the value of $\delta_3$.
Rules for evaluating $\W(\infty)$ are detailed in \ref{sec:counting-fpps}.
\begin{enumerate} [{Case} 1{:}, leftmargin=0.5in]
    \item $\delta_3 = 0$.\\
        By Lemma \ref{lemma:delta3}, this case occurs when         $\W(\infty)=\Y(\infty)=\infty$, that is, when $v>u$ and $r>s$.
        In particular, for this condition, $\delta_1=d$ exactly when
        $v=r=(d+1)/2$ and $u=s=(d-1)/2$, and $\delta_1<d$ otherwise.
        Thus $\delta_1+\delta_3\le d$.
    \item $\delta_3 = 1$.\\
        By Lemma \ref{lemma:delta3}, this case occurs either when  $\W(\infty)=\infty$ and $\Y(\infty)=a \in \Fq$ (that is, $v>u$ and $r \leq s$), or similarly, when $\Y(\infty)=\infty$ and $\W(\infty)=a \in \Fq$, 
        (that is, $r>s$ and $v \leq u$).
        Either way,  $\delta_1 \leq d-1$, 
        so $\delta_1+\delta_3\le d$.
    \item $\delta_3 = 2$.\\
        By Lemma \ref{lemma:delta3}, this case occurs when $\W(\infty)=a$ for some $a \in \Fq $  and $\Y(\infty)=b$ for some $b\in \Fq$ (that is, when $v \leq u$ and $r \leq s$).
        This means that $\delta_1 \leq (d-3)/2+(d-3)/2$, 
        so $\delta_1+\delta_3\le d$.
\end{enumerate} \end{proof}


\begin{thm}
\label{thm:M(q,q-d)}
For odd $d \ge 3$ and $q \ge d+2,~M(q,q-d) \ge \S_d(q).$
\end{thm}
\begin{proof}
This follows from Lemma \ref{lemma:hd(Pi-contraction)}, because by definition, $\S_d(q) = |\bigcup_{v,u} \Pi'_{v,u}(q)|$, for $v$ and $u$ as specified in Definition \ref{dfn:S_d(q)} and in Lemma \ref{lemma:hd(Pi-contraction)}.
\end{proof}

\bigskip
We have derived explicit formulas for $\S_5(q)$ and $\S_7(q)$ based the  congruence classes $q \Mod 3$, by substituting known formulas for each term in the definitions of $\S_5(q)$ and $\S_7(q)$. 
The derived formulas for $\S_5(q)$ and $\S_7(q)$ are 
given in Theorems \ref{thm:M(q,q-5)formulas} and \ref{thm:M(q,q-7)formulas} below.

\begin{thm}
\label{thm:M(q,q-5)formulas}
For $q>7$
    \begin{equation*}
    M(q,q-5) \ge \S_5(q)=
    \begin{cases}
        \frac{1}{2}(q^4+q^3+q^2-3q),& \mathrm{~for~} q \equiv{0 \Mod 3},\\
        \frac{1}{2}(q^4+q^2-2q), & \mathrm{~for~odd~} q \equiv{1\Mod 3}, \\
        \frac{1}{2}(q^4+3q^2-4q), & \mathrm{~for~even~} q \equiv{1 \Mod 3},\\
        \frac{1}{2} (q^4+2q^3-q^2-2q), & \mathrm{~for~odd~} q \equiv{2 \Mod 3},\\
        \frac{1}{2} (q^4+2q^3+q^2-4q), & \mathrm{~for~even~}q \equiv{2 \Mod 3}.
    \end{cases}
\end{equation*}
\end{thm}

\begin{proof}
We use Theorem \ref{thm:M(q,q-d)}.
We provide a detailed example of the derivation of the formula for $\S_5(q)$ when $q$ is odd and  $q \equiv{2 \Mod 3}$.

\begin{align*}
    \S_5(q) 
    &= \S_3(q) +  \sum_{u<3}N_{3,u}(q) + \sum_{v \le 1}{N_{v,1}(q) }& \mathrm{(~Equation~15)}\\
    &=\S_3(q) + (N_{3,2}(q)+N_{3,0}(q)) + (N_{1,1}(q)+N_{0,1}(q))&\mathrm{(Observation~4)}\\
    &=q(q-1) + \tfrac{1}{2}q^2(q-1)^2 + q^2(q-1) + q(q-1)^2 + q(q-1)\\
    &&\mathrm{(Equations~2,~4,~7,~10~and~Theorem~2.3)}\\
    &=\tfrac{1}{2}(q^4+2q^3-q^2-2q)
\end{align*}
The other results are obtained in a similar way by adding terms appropriate for each congruence class. \end{proof}

\begin{thm}
\label{thm:M(q,q-7)formulas}
For $q>9$,
    \begin{equation*}
    M(q,q-7) \ge \S_7(q) = 
    \begin{cases}
        \frac{1}{6}(2q^5+q^4+q^3+5q^2-9q),
        & \mathrm{~for~} q \equiv{0 \Mod 3},\\
        \frac{1}{6}(2q^5+q^4-2q^3+5q^2-6q), 
        & \mathrm{~for~odd~} q \equiv{1\Mod 3}, \\
        \frac{1}{6}(2q^5+3q^4+2q^3+9q^2-16q), 
        & \mathrm{~for~even~} q \equiv{1 \Mod 3},\\
        \frac{1}{6} (2q^5+q^4+4q^3-q^2-6q), 
        & \mathrm{~for~odd~} q \equiv{2 \Mod 3},\\
        \frac{1}{6} (2q^5+3q^4+8q^3+3q^2-16q), 
        & \mathrm{~for~even~}q \equiv{2 \Mod 3}.
    \end{cases}
\end{equation*}
\end{thm}

\begin{proof}
We use Theorem \ref{thm:M(q,q-d)}.
We provide a detailed example of the derivation of the formula for $\S_7(q)$ when $q$ is even and  $q \equiv{2 \Mod 3}$.

\begin{align*}
    \S_7(q) 
    &= \S_5(q) +  \sum_{u<4}N_{4,u}(q) + \sum_{v \le 2}{N_{v,2}(q) } ~~~~~~~~~~~~~~~~~~~~~~~~~~~~~~~\mathrm{(E quation~16)}\\
    &=\S_5(q) + (N_{4,3}(q)+N_{4,0}(q)) +
    (N_{2,2}(q)+N_{2,0}(q))
    ~~~~~~~~~~\mathrm{(Observation~4~and~Equation~9)}\\
    &=\frac{1}{2}(q^4+2q^3+q^2-4q) + \frac{1}{3} (q+1)q^2(q-1)^2 + \frac{1}{3}q(q-1)(q^2+2) + q(q-1)^2 + q(q-1)\\
    &~~~~~~~~~~~~~~~~~~~~~~~~~~~~~~~~~~~~~~~~~~~~~~~~~~~~~~~~~~~~~~~~~~~~\mathrm{(Th.~2.3,~2.5~and~4.6,~ and~Eq.~2~and~5)}\\
    &=\tfrac{1}{6} (2q^5+3q^4+8q^3+3q^2-16q).
\end{align*}

Derivation of formulas for the other cases is similar, using terms appropriate for each congruence class.
\end{proof}




To illustrate, for $q=59$ and $d=5$, Theorem \ref{thm:M(q,q-5)formulas} gives $M(59,54) \ge
6,262,260$, which is an improved lower bound. The previous lower bound of
6,067,206 was obtained from \fs for $q=59$ \cite{bmmssPRFwaifi2020}.
For $q = 59$ and $d = 7$, Theorem \ref{thm:M(q,q-7)formulas} gives $M(59,52) \ge 240,463,940$, which is an improved lower bound, where the previous lower bound was 240,262,042, using \fs \cite{bmmssPRFwaifi2020}.

Theorems \ref{thm:M(q,q-d)}, \ref{thm:M(q,q-5)formulas} and \ref{thm:M(q,q-7)formulas} give new lower bounds on $M(q,q-d)$ for odd $d$.
So again, we relax the Hamming distance constraint in order to consider an additional set of \fs for results for $M(q,q-d-1)$, \ie{when the desired Hamming distance is $q-d-1$ and, hence, $d-1$ is even}.
We investigate the effect of contraction on the Hamming distance of the set of permutations $\Pi^m_{\frac{d+3}{2},\frac{d+1}{2}}(q)$.
These are the permutations produced by the monic \fs, \ie{the set $P^m_{\frac{d+3}{2},\frac{d+1}{2}}(q)$}. 
By Definition \ref{dfn:counting}\ref{dfn:N^m_{vu}}, there are $N^m_{\frac{d+3}{2},\frac{d+1}{2}}(q)$ permutations in this set.
Let  $v,u,s \le \frac{d+1}{2}$ and $r=\frac{d+3}{2}$.
Lemma \ref{lemma:hd(monicPi-contraction)} shows that $\Pi^{'m}_{r,s}(q)$, the set of permutations  generated by the monic \fs 
$P'^m_{r,s}(q)$ 
has Hamming distance $q-d-1$.
Lemma \ref{lemma:hd(Pi-contraction,monicPi-contraction)} shows that permutations in the union of the two sets $\bigcup_{v,u} \Pi'_{v,u}(q)$  and $\Pi^{'m}_{r,s}(q)$ have pairwise Hamming distance at least $q-d-1$.
Theorem \ref{thm:M(q,q-d-1)} uses these results to give the new lower bound $M(q,q-d-1) \ge \S_d(q)+N^m_{\frac{d+3}{2},\frac{d+1}{2}}(q)$, for $q>7$ and for $d \ge 3$.
We give explicit formulas for $M(q,q-6)$ in Theorem \ref{thm:M(q,q-6)formulas}.

\begin{lemma}
\label{lemma:hd(monicPi-contraction)}
For odd $d \ge 3$, let $r=(d+3)/2$ and $s=(d+1)/2$.
Let $\Pi^{'m}_{r,s}(q)$ be the set of permutations generated by the operation of contraction on $\Pi^m_{r,s}(q))$.
Then $hd(\Pi^{'m}_{r,s}(q)) \ge q-d-1.$
\end{lemma}

\begin{proof}
By Lemma \ref{lemma:hd-contraction}, it suffices to show that $\delta_1+\delta_3\le d+1$.
Note that $\delta_1=\deg(V(x)S(x) - U(x)R(x))\le d+2$. 
However, since the \fs $P^m_{r,s}(q)$ are monic, the high order term of the polynomial
$(V(x)S(x) - U(x)R(x))$ is zero, so the in fact, $\delta_1=\deg(V(x)S(x) - U(x)R(x))\le d+1$.
Note also that since $r>s, ~\W(\infty)=\Y(\infty)=\infty$.
Thus $\delta_3=0$, and $\delta_1+\delta_3\le d+1$.
\end{proof}

\begin{lemma}
\label{lemma:hd(Pi-contraction,monicPi-contraction)}
For odd $d \ge 3$, let $r=(d+3)/2$~and $s=(d+1)/2$. Let $v$ and $u$ be defined by 
  \begin{equation*}
    v,u=
    \begin{cases}
        v \le (d+1)/2,
        &\mathrm{~when~} v>u, \\
 
        u \le (d-3)/2, & \mathrm{~when~} v<u, \mathrm{~and}\\
        u,v \le (d-3)/2,& \mathrm{~when~} v=u.
    \end{cases}
  \end{equation*} 
Then $hd(\bigcup_{v,u} \Pi'_{v,u}(q),\Pi^{'m}_{r,s}(q)) \ge q-d-1.$
\end{lemma}

\begin{proof}
Let $\pi \in \bigcup_{v,u} \Pi_{v,u}(q)$ and $\sigma \in \Pi^m(v,u)(q)$.
Let $\W(x)=\frac{V(x)}{U(x)}$ and $\Y(x)=\frac{R(x)}{S(x)}$ be the \fs that generate $\pi$ and $\sigma$, respectively. Note that $r$ and $s$ are fixed at $r=(d+3)/2$ and $s=d+1)/2$.
We show that $hd(\pi',\sigma') \geq q-d-1.$
By Lemma \ref{lemma:hd-contraction}, it suffices to show that $\delta_1+\delta_3\le d+1$.
We do a proof by cases based on the value of $\delta_3$.
\begin{enumerate} [{Case} 1{:}, leftmargin=0.5in]
    \item $\delta_3 = 0$.\\
        By Lemma \ref{lemma:delta3}, this case occurs when  $\W(\infty)=\Y(\infty)=\infty$, that is, when $v>u$ and $r>s$.
        In particular, $\delta_1=d+1$ exactly when when $u,v,r$ and $s$ achieve their maximum values), \ie{
        $v=(d+1)/2$, $u=(d-1)/2$, $r=(d+3)/2$ and $s=(d+1)/2$}.
        $\delta_1<d+1$ otherwise.
        Thus $\delta_1+\delta_3\le d+1$.
        
    \item $\delta_3 = 1$.\\
        By Lemma \ref{lemma:delta3}, this case occurs when $r>s$ and $v \leq u$, that is, when $\Y(\infty)=\infty$ and $\W(\infty)=a \in \Fq$. 
        So $\delta_1=\max((d-3)/2+(d+1)/2,(d-3)/2+(d+3)/2)=d$.
        Thus $\delta_1+\delta_3\le d+1$.
 
    \item $\delta_3 = 2$.\\
        By Lemma \ref{lemma:delta3}, this case occurs when  $\W(\infty)=a$ for some $a \in \Fq $  and $\Y(\infty)=b$ for some $b\in \Fq$ (that is, when $v \leq u$ and $r \leq s$).
        This means that $\delta_1 \leq ((d-3)/2+(d+1)/2,(d+1)/2+(d-1)/2)=d-1$, 
        so $\delta_1+\delta_3 \le d+1$.

\end{enumerate}
\end{proof}


\begin{thm}
\label{thm:M(q,q-d-1)}
For $d \ge 3$ and $q \ge d+2$, $M(q,q-d-1) \ge \S_d(q)+N^m_{\frac{d+3}{2},\frac{d+1}{2}}(q)$
\end{thm}

\begin{proof}
This follows from Theorems \ref{lemma:hd(Pi-contraction)}, \ref{lemma:hd(monicPi-contraction)} and \ref{lemma:hd(Pi-contraction,monicPi-contraction)}, since $\S_d(q)=|\bigcup_{v,u} \Pi'_{v,u}(q)|$, and $N^{'m}_{r,s}(q)=|\Pi^m_{r,s}(q)|$, for $u,v,s$ and $r$ as given.
\end{proof}

\begin{thm}
\label{thm:M(q,q-6)formulas}
For $q>9$
    \begin{equation*}
    M(q,q-6) \ge \S_5(q) + N^m_{4,3}(q) = 
    \begin{cases} 
        \frac{1}{6}(5q^4+3q^3+q^2-9q),
        & \mathrm{~for~} q \equiv{0 \Mod 3},\\
        \frac{1}{6}(5q^4+q^2-6q), 
        & \mathrm{~for~odd~} q \equiv{1\Mod 3}, \\
        \frac{1}{6}(5q^4+7q^2-12q), 
        & \mathrm{~for~even~} q \equiv{1 \Mod 3},\\
        \frac{1}{6} (5q^4+6q^3-5q^2-6q), 
        & \mathrm{~for~odd~} q \equiv{2 \Mod 3},\\
        \frac{1}{6} (5q^4+6q^3+q^2-12q), 
        & \mathrm{~for~even~}q \equiv{2 \Mod 3}.
    \end{cases}
\end{equation*}
\end{thm}

\begin{proof}
We use Theorem \ref{thm:M(q,q-d-1)}.
As mentioned earlier, the operation of contraction on a PA does not affect the size of the PA, that is, $|\Pi^{'m}_{r,s}(q)|=|\Pi^{m}_{r,s}(q)|=N^m_{4,3}(q)$.
We provide a detailed example of the derivation of the formula for $\S_5(q)+N^m_{4,3}(q)$ when $q$ is odd and  $q \equiv{2 \Mod 3}$.

\begin{align*}
    M(q,q-6) &\ge \S_5(q) + N^m_{4,3}(q)&\mathrm{(Theorem~4.10)}\\
    & \ge \S_5(q)+ \tfrac{1}{q-1}N_{4,3}(q)&
    \mathrm{(Observation~3)}\\
    &=  \tfrac{1}{2} (q^4+2q^3-q^2-2q) + \tfrac{1}{q-1} ((q+1)q^2(q-1)^2)/3&\mathrm{(Theorems~4.6~and~2.5)} \\
    &=\tfrac{1}{6} (5q^4+6q^3-5q^2-6q).
\end{align*}

\noindent The other formulas are obtained in a similar way by adding terms appropriate for each congruence class. \end{proof}


For example, for $q=59$, Theorem \ref{thm:M(q,q-6)formulas} gives the improved lower bound $M(59,53) \ge 10,300,220$. The previous lower bound was $M(59,53) \ge 407,218$, given by permutation polynomials \cite{bmmssw-19PermPolynomials}.




\section{Lower bounds for $M(q,q-D)$ and $M(q+1,q-D),~D \ge 8$. }
\label{sec:conjectures}
A simple computational strategy computes $N_{v,u}(q)$ by considering all possible rational functions, $\frac{V(x)}{U(x)}$, where $V(x)$ and $U(x)$ are polynomials of degree $v$ and $u$, respectively. 
This entails evaluating $q^{u+v+2}$ different rational functions. 
Using normalization results and the $F$ and  $G$ maps \cite{bmmssPRFwaifi2020}, five coefficients can be fixed, thus reducing the search to at most $q^{u+v-3}$ different rational functions, an efficiency that allows additional computational results. 
 
For instance, we have computed $N_{5,4}(q)$ for various values of $q$, which we use as a basis for the following theorem.

\begin{thm} \label{thm:N_{5,4}}
For $32 \le q \le 127$, 
\begin{equation*}
    N_{5,4}(q) =
   \begin{cases}
        \frac{1}{2}(q^6-2q^4+62q^3-61q^2)& \mathrm{~when~} q \equiv{0 \Mod 5},\\
        \frac{1}{2}(q^6-q^4-2q^3+2q^2), & \mathrm{~when~} q \equiv{1,4 \Mod 5},  \mathrm{~and}\\
        \frac{1}{2}(q^6-2q^4+q^2),
        &\mathrm{~when~} q \equiv{2,3 \Mod 5}.
    \end{cases}
\end{equation*} \qed
\end{thm} 

We also computed $N_{5,4}(q)$, for $q \le 31$.
The results are listed in Table \ref{tbl:N5,4andN5,5andT9}. $N_{5,4}(q)$ entries marked with an asterisk denote agreement with Theorem \ref{thm:N_{5,4}}.
We believe the formulas given in Theorem \ref{thm:N_{5,4}} are also valid for all $q
\ge 128$. This is given in Conjecture \ref{5,4 conj}.

\begin{conjecture}
\label{5,4 conj}
Theorem \ref{thm:N_{5,4}} is true, for all $q \ge 128$.
\end{conjecture}

To compute $N_{5,5}(q)$, we use Theorem \ref{thm:equal-num-denom} and Theorem \ref{thm:N_{5,4}} as follows. 
First, Observation \ref{obs:N_v1=0} gives $N_{5,1}(q)=0$. 
In addition, we observe experimentally that  $N_{5,2}(q)=0$ and $N_{5,3}(q)=0$,  for all $32 \le q
\le 127$. 
So by Theorem \ref{thm:equal-num-denom}, $N_{5,5}(q)$ = $(q-1)(N_{5,4}(q)+N_{5,0}(q))$ for $q$ within this range. 
Using the formulas for $N_{5,0}(q)$ given  by Equation \ref{eqn:N_5,0},
we obtain the explicit formulas for $N_{5,5}(q)$ in Theorem \ref{thm:N_{5,5}}. 
The derivation of these formulas is similar to the derivations described for other formulas given earlier in this paper. 
(See for example Theorems \ref{thm:M(q+1,q-4)formulas}, \ref{thm:M(q,q-5)formulas}, etc.).

\bigskip

\begin{table}[h!tb]
\centering
\renewcommand{\arraystretch}{1.4}
\begin{tabular}{|c|c|c|c|}
\hline $\boldsymbol{q}$ & $\boldsymbol{N_{5,4}(q)}$ & $\boldsymbol{N_{5,5}(q)}$&$\boldsymbol{\T_9(q)}$   \\
\hline 

 16& 12,110,400 & 182,808,000&213,579,840\\
\hline
17&16,189,44&266,207,488&308,227,680\\
\hline
19&24,503,958 & 446,802,480 &513,218,880\\
\hline 
23&74,762,512  & 1,650,664,092&1,853,247,264 \\
\hline
25&125,820,000 & 3,030,033,600 &3,368,352,000\\
\hline
27&193,179,168$^*$ & 5,035,964,076 &5,559,326,136\\
\hline
29&297,858,652$^*$ & 8,340,701,600 &9,144,880,160\\
\hline
31&444,126,150 & 13,332,433,500 &14,530,349,760\\
\hline
\end{tabular}
\caption{$N_{5,4}(q)$, $N_{5,5}(q)$, and $\T_9(q)$ for $q \le 31$, obtained by experimental computation. ($^*$ denotes agreement with $N_{5,4}$ conjecture.)}
\label{tbl:N5,4andN5,5andT9}
\end{table}


\begin{thm}
\label{thm:N_{5,5}}
For $32 \le q \le 127$, 
\begin{equation*}
N_{5,5}(q) =
\begin{cases}
        \frac{1}{4}(2q^7-2q^6-2q^5+125q^4-245q^3+121q^2+q),
        & \mathrm{~when~} q \equiv{0\Mod 5}\\
        \frac{1}{2}(q^7-q^6-q^5-q^4+4q^3-2q^3-2q^2), 
        & \mathrm{~when~} q \equiv{1 \Mod 5},  \\
        \frac{1}{2}(q^7-q^6-2q^4+3q^3-q^2),
        &\mathrm{~when~} q \equiv{2,3 \Mod 5}  , \\
        \frac{1}{2}(q^7-q^6-q^5+q^4),
        & \mathrm{~when~} q \equiv{4 \Mod 5}
    \end{cases}
\end{equation*}  \qed
\end{thm}

\begin{conj}
Theorem \ref{thm:N_{5,5}} is true, for all $q \ge 128$.
\label{conj:N_5,5-conj}

\end{conj}

It should be noted this is valid if Conjecture \ref{5,4 conj} is valid.

\begin{thm} 
\label{thm:T_9 theorem}
For $32 \le q \le 127, M(q+1,q-9) \ge \T_9(q)$, where $\T_9(q)$ is given by \\

For odd $q\equiv{0\Mod 3}$:
\begin{equation*} 
\T_9(q)=
\begin{cases}                   
        \frac{1}{6}(3q^7+5q^6-13q^4+6q^3+8q^2-9q), 
        & q\equiv{1 \Mod 5},\\
        \frac{1}{6}(3q^7+5q^6+3q^5-10q^4+3q^3+5q^2-9q),
        & q\equiv{2,3 \Mod 5},\\
        \frac{1}{6}(3q^7+5q^6-7q^4+6q^3+2q^2-9q), 
        & q\equiv{4 \Mod 5},\\
    \end{cases}
\end{equation*}

For odd $q\equiv{1\Mod 3}$:
\begin{equation*} 
\T_9(q)=
    \begin{cases}
        \frac{1}{12}(6q^7+10q^6+349q^4+9q^3-359q^2-15q),  
        & q\equiv{0 \Mod 5},\\
        \frac{1}{6}(3q^7+5q^6-16q^4+3q^3+11q^2-6q), & q\equiv{1\Mod 5},\\
        \frac{1}{6}(3q^7+5q^6+3q^5-13q^4+8q^2-6q),
        &q\equiv{2,3\Mod 5},\\
        \frac{1}{6}(3q^7+5q^6-10q^4+3q^3+5q^2-6q), & q\equiv{4 \Mod 5}, \\
    \end{cases}
\end{equation*}

For odd $q\equiv{2\Mod 3}$:
\begin{equation*} 
\T_9(q)=
    \begin{cases}
        \frac{1}{12}(6q^7+10q^6+361q^4+9q^3-371q^2-15q), 
        & q\equiv{0 \Mod 5},\\
        \frac{1}{6}(3q^7+5q^6-10q^4+3q^3+5q^2-6q),
        &q\equiv{1 \Mod 5},\\
        \frac{1}{6}(3q^7+5q^6+3q^5-7q^4+2q^2-6q),
        & q\equiv{2,3 \Mod 5},\\
        \frac{1}{6}(3q^7+5q^6-4q^4+3q^3-q^2-6q),
        & q\equiv{4 \Mod 5},\\
    \end{cases}
\end{equation*}

For even $q\equiv{1\Mod 3}$:
\begin{align*} 
    \T_9(q)= \tfrac{1}{6}(3q^7+5q^6+2q^5-10q^4+11q^3+5q^2-16q),
     ~~~ q\equiv{1,4 \Mod 5}, 
\end{align*}

For even $q\equiv{2\Mod 3}$:
\begin{align*}
    \T_9(q)= \tfrac{1}{6}(3q^7+5q^6+5q^5-7q^4+8q^3+2q^2-16q),
    ~~~~~~~q\equiv{2,3 \Mod 5}.
\end{align*} 

Note: Empty congruence classes are not listed.\qed
\label{thm:M(q+1,q-9)}
\end{thm} 


It should be noted that the equations indicated in Theorems \ref{thm:T_9 theorem} have been verified by computation. For example, for $q=59$, Theorem \ref{thm:M(q+1,q-9)} gives the improved lower bound $M(60,50) \ge 1,279,468,210,920$. The previous lower bound was $M(60,50) \ge 363,621,720$, given by \fs \cite{bmmssPRFwaifi2020}.

\begin{conj}
Theorem \ref{thm:M(q+1,q-9)} is true, for all $q \ge 128$.
\label{conj:M(q+1,q-9)}

\end{conj}

It should be noted this is valid if Conjecture \ref{5,4 conj} is valid.

\begin{thm} 
\label{thm:M(q,q-9)}
For $32 \le q \le 127$, $M(q,q-9) \ge \S_9(q)$, where $\S_9(q)$ is given by \\

For odd $q\equiv{0\Mod 3}$:
\begin{equation*} 
\S_9(q)=
\begin{cases}                   
        \frac{1}{6}( 3q^6 + 5q^5 -  5q^4 - 2q^3 + 8q^2 - 9q),
        & q\equiv{1\Mod 5},\\
        \frac{1}{6}( 3q^6 + 5q^5 - 2q^4 + q^3 - q^2 - 6q),
        & q\equiv{2,3 \Mod 5},\\
        \frac{1}{6}(3q^6 + 5q^5 - 8q^4 + 13q^3 - 7q^2 - 6q),
        & q\equiv{4 \Mod 5},\\
    \end{cases}
\end{equation*}

For odd $q\equiv{1\Mod 3}$:
\begin{equation*} 
\S_9(q)=
    \begin{cases}
        \frac{1}{6}( 3q^6 + 5q^5 -  8q^4 - 5q^3 + 11q^2 - 6q), 
        & q\equiv{1\Mod 5},\\
        \frac{1}{6}(3q^6 + 5q^5 -  5q^4 - 5q^3 + 8q^2 - 6q),
        &q\equiv{2,3\Mod 5},\\
        \frac{1}{6}(3q^6+5q^5-8q^4+q^3+5q^2-6q), & q\equiv{4 \Mod 5}, \\
    \end{cases}
\end{equation*}

For odd $q\equiv{2\Mod 3}$:
\begin{equation*} 
\S_9(q)=
    \begin{cases}
        \frac{1}{12}( 6q^6 + 10q^5 - 4q^4 + 371q^3 - 368q^2 - 15q), 
        & q\equiv{0 \Mod 5},\\
        \frac{1}{6}( 3q^6 + 5 q^5 - 2q^4 - 5 q^3 + 5 q^2 - 6q),
        &q\equiv{1 \Mod 5},\\
        \frac{1}{6}( 3q^6 + 5 q^5 + q^4- 5q^3 + 2q^2 - 6q),
        & q\equiv{2,3 \Mod 5},\\
        \frac{1}{6}( 3q^6 + 5q^5 - 2q^4 + q^3 - q^2 - 6q),
        & q\equiv{4 \Mod 5},\\
    \end{cases}
\end{equation*}

For even $q\equiv{1\Mod 3}$:
\begin{align*} 
\S_9(q)=     
    \begin{cases}
        \frac{1}{6}(3q^6 + 5q^5 - 6q^4 - q^3 + 15q^2 - 16q),
        & q\equiv{1\Mod 5},\\
        \frac{1}{6}(3q^6 + 5q^5 - 6q^4 +  5q^3 + 9q^2 - 16q),
        & q\equiv{4 \Mod 5}.
    \end{cases}
\end{align*}

For even $q\equiv{2\Mod 3}$:
\begin{align*}
    \S_9(q)= \tfrac{1}{6}(3q^6 + 5q^5 + 3q^4 - q^3 + 6q^2 - 16q),
    ~~~~~~~q\equiv{2,3 \Mod 5}.
\end{align*}

Note: Empty congruence classes are not listed. \qed
\end{thm} 


It should be noted that the equations indicated in Theorem \ref{thm:M(q,q-9)} have been verified by computation.
For example, for $q=59$, Theorem  \ref{thm:M(q,q-9)} gives the improved lower bound $M(59,50) \ge 21,682,031,540$. The previous lower bound was $M(59,50) \ge 20,979,628,398$, given by \fs \cite{bmmssPRFwaifi2020}.

\begin{conj}
Theorem \ref{thm:M(q,q-9)} is true for all $q \ge 128$.
\label{conj:M(q,q-9)}
\end{conj}

It should be noted this is valid if Conjecture \ref{5,4 conj} is valid.



Similarly, we obtained the following theorems and conjectures by computation:

\begin{thm} 
\label{thm:M(q+1,q-8)}
For $32 \le q \le 127$, $M(q+1,q-8) \ge \T_7(q)+N^m_{5,4}(q)$, where $\T_7(q)+N^m_{5,4}(q)$ is given by \\

For odd $q\equiv{0\Mod 3}$:
\begin{equation*} 
\T_7(q)+N^m_{5,4}(q)=
\begin{cases}                   
        \frac{1}{6}(2q^6+6q^5-q^4+6q^3-4q^2-9q) 
        & q\equiv{1,4 \Mod 5},\\
        \frac{1}{6}(2q^6+6q^5-q^4+3q^3-q^2-9q)
        & q\equiv{2,3 \Mod 5},\\
        
    \end{cases}
\end{equation*}

For odd $q\equiv{1\Mod 3}$:
\begin{equation*} 
\T_7(q)+N^m_{5,4}(q)=
    \begin{cases}
        \frac{1}{3}(q^6+3q^5-2q^4+94q^2-3q),  
        & q\equiv{0 \Mod 5},\\
        \frac{1}{6}(2q^6+6q^5-4q^4+3q^3-q^2-6q), & q\equiv{1,4\Mod 5},\\
        \frac{1}{3}(q^6+3q^5-2q^4+q^2-3q),
        &q\equiv{2,3\Mod 5},\\
    \end{cases}
\end{equation*}

For odd $q\equiv{2\Mod 3}$:
\begin{equation*} 
\T_7(q)+N^m_{5,4}(q)=
    \begin{cases}
        \frac{1}{3}(q^6+3q^5+q^4+91q^2-3q), 
        & q\equiv{0 \Mod 5},\\
        \frac{1}{6}(2q^6+6q^5+2q^4+6q^3-7q^2-6q),
        &q\equiv{1,4 \Mod 5},\\
        \frac{1}{3}(q^6+3q^5+q^4-2q^2-3q),
        & q\equiv{2,3 \Mod 5},\\
    \end{cases}
\end{equation*}

For even $q\equiv{1\Mod 3}$:
\begin{align*} 
    \T_7(q)+N^m_{5,4}(q)= \tfrac{1}{6}(2q^6+8q^5-4q^4+11q^3-q^2-16q),
     ~~~ q\equiv{1,4 \Mod 5}, 
\end{align*}

For even $q\equiv{2\Mod 3}$:
\begin{align*}
    \T_7(q)+N^m_{5,4}(q)= \tfrac{1}{3}(q^6+4q^5+q^4+4q^3-2q^2-8q),
    ~~~~~~~q\equiv{2,3 \Mod 5}.
\end{align*} 

Note: Empty congruence classes are not listed. \qed
\end{thm}

For example, for $q=47$, Theorem  \ref{thm:M(q+1,q-8)} gives the improved lower bound $M(47,39) \ge
1,94,389,744$. The previous lower bound was $M(47,39) \ge 77,228,802$, given by \fs \cite{bmmssPRFwaifi2020}.

\begin{conj}
Theorem \ref{thm:M(q+1,q-8)} is true, for all $q \ge 128$. \qed
\label{conj:M(q+1,q-8)}
\end{conj}


\begin{thm} 
\label{thm:M(q,q-8)}
For 
$32 \le q \le 127$, $M(q,q-8) \ge \S_7(q)+N^m_{5,4}(q)$,\\ where $\S_7(q)+N^m_{5,4}(q)$ is given by \\

For odd $q\equiv{0\Mod 3}$:
\begin{equation*} 
\S_7(q)+N^m_{5,4}(q)=
\begin{cases}                   
        \frac{1}{6}(5q^5+4q^4+q^3-q^2-9q) 
        & q\equiv{1,4 \Mod 5},\\
        \frac{1}{6}(5q^5+4q^4-2q^3+2q^2-9q)
        & q\equiv{2,3 \Mod 5},\\
        
    \end{cases}
\end{equation*}

For odd $q\equiv{1\Mod 3}$:
\begin{equation*} 
\S_7(q)+N^m_{5,4}(q)=
    \begin{cases}
        \frac{1}{6}(5q^5+4q^4-5q^3+188q^2-6q),  
        & q\equiv{0 \Mod 5},\\
        \frac{1}{6}(5q^5+4q^4-2q^3-q^2-6q), & q\equiv{1,4\Mod 5},\\
        \frac{1}{6}(5q^5+4q^4-5q^3+2q^2-6q),
        &q\equiv{2,3\Mod 5},\\
    \end{cases}
\end{equation*}

For odd $q\equiv{2\Mod 3}$:
\begin{equation*} 
\S_7(q)+N^m_{5,4}(q)=
    \begin{cases}
        \frac{1}{12}(5q^5+4q^4+q^3+182q^2-6q), 
        & q\equiv{0 \Mod 5},\\
        \frac{1}{6}(5q^5+4q^4+4q^3-7q^2-6q),
        &q\equiv{1,4 \Mod 5},\\
        \frac{1}{6}(5q^5+4q^4+q^3-4q^2-6q),
        & q\equiv{2,3 \Mod 5},\\
    \end{cases}
\end{equation*}

For even $q\equiv{1\Mod 3}$:
\begin{align*} 
    \S_7(q)+N^m_{5,4}(q)= \tfrac{1}{6}(5q^5+6q^4+2q^3+3q^2-16q),
     ~~~ q\equiv{1,4 \Mod 5}, 
\end{align*}

For even $q\equiv{2\Mod 3}$:
\begin{align*}
    \S_7(q)+N^m_{5,4}(q)= \tfrac{1}{6}(5q^5+6q^4+5q^3-16q),
    ~~~~~~~q\equiv{2,3 \Mod 5}.
\end{align*} 

Note: Empty congruence classes are not listed. \qed
\end{thm} 

For example, for $q=47$, Theorem  \ref{thm:M(q,q-8)} gives the improved lower bound $M(47,39) \ge
1,94,389,744$. The previous lower bound was $M(47,39) \ge 77,228,802$, given by \fs \cite{bmmssPRFwaifi2020}.

\begin{conj}
Theorem \ref{thm:M(q,q-8)} is true, for all $q \ge 128$. 
\label{conj:M(q,q-8)}
\end{conj}

\bigskip
We have also computed 
$N_{6,5}(q)$, $N_{6,4}(q)$, $N_{6,3}(q)$, and $N_{6,2}(q)$, for prime powers $q,~ 16 \le q \le 41$. These are shown in Table \ref{tbl:N6,5}. We also found
$N_{7,6}(19)=19,960,510,914$.

\begin{table}[!htb]
\centering
\renewcommand{\arraystretch}{1.4}
\begin{tabular}{|c|c|c|c|c|}
\hline $\boldsymbol{q}$ & $\boldsymbol{N_{6,5}(q)}$&$\boldsymbol{N_{6,4}(q)}$& $\boldsymbol{N_{6,3}(q)}$&$\boldsymbol{N_{6,2}(q)}$\\
\hline 

16& 204,963,840&17,524,800&460,800& 86,400\\
\hline
17&84,489,728&5,548,800&443,904&73,984\\
\hline
19&33,334,740&2,573,208&116,964&0\\
\hline 
23&0 &0&0&0\\
\hline
25&1,512,000&360,000&0&0\\
\hline
27&12,320,100&492,804&0&492,804\\
\hline
29&3,956,064&0&0&0\\
\hline
31&0&0&0&0\\
\hline
32&291,282,944&24,847,616&0&0\\
\hline
37&0&0&0&0\\
\hline
41&0&0&0&0\\
\hline

\end{tabular}
\caption{$N_{6,5}(q)$, $N_{6,4}(q)$, $N_{6,3}(q)$, and $N_{6,2}(q)$, for $ 16 \le q \le 41$.}
\label{tbl:N6,5}
\end{table}

\section{Conclusions and Final Conjectures}
\label{sec:conclusions}

Using recent enumerations of all degree 3 \fs, and formulas for the number of all degree 3 \fs, given by
Ferraguti and Micheli \cite{Ferraguti20}, and for degree 4 \fs by Ding and Zieve \cite{Ding2020}, and separately by Hou \cite{Hou20}, we have given improved lower bounds for $M(q,q-d)$ and
$M(q+1,q-d)$, for $q>9$, and for $d \in \{4,5,6,7\}$ 

Using computation we described formulas for $N_{5,4}(q)$, for all $32 \le q \le 127$. We conjecture that the formulas are valid also for all $q \ge 128$.
Using the formulas for $N_{5,4}(q)$ we gave a formula for $N_{5,5}(q)$ and formulas for improved lower bounds for $M(q,q-d)$ and $M(q+1,q-d)$, for all $32 \le q \le 127$, and for $d \in \{8,9\}$. Again, we conjecture that these formulas are also valid for all $q \ge 128$. 

Based on our computations, we conjecture that $N_{6,u}(q) = 0$, for all $q \ge 37$ and all $1 \le u \le 5$.
Computed values for
$N_{6,u}(q) = 0$ are given in Table \ref{tbl:N6,5}.

\bibliographystyle{plainnat} 

\end{document}